\magnification=1095 \baselineskip=15pt

    \hoffset=-17mm 

\def\today{\number\day\space\ifcase\month\or January\or February\or March\or
  April\or May\or June\or July\or August\or September\or October\or
  November\or December\fi\space\number\year}


\font\ninerm=cmr9   \font\eightrm=cmr8    \font\sixrm=cmr6
\font\ninei=cmmi9   \font\eighti=cmmi8    \font\sixi=cmmi6
\font\ninesy=cmsy9   \font\eightsy=cmsy8    \font\sixsy=cmsy6
\font\ninebf=cmbx9   \font\eightbf=cmbx8    \font\sixbf=cmbx6
\font\ninett=cmtt9  \font\eighttt=cmtt8
\font\nineit=cmmi9  \font\eightit=cmmi8
\font\ninesl=cmsl9  \font\eightsl=cmsl8
\font\tenbcal=cmbsy10
\font\ninebcal=cmbsy10 at9pt
\font\eightbcal=cmbsy10 at8pt
\font\sevenbcal=cmbsy10 at7pt
\font\sixbcal=cmbsy10 at6pt
\font\fivebcal=cmbsy10 at5pt
\newfam\bcalfam
\def\bcal#1{{\fam\bcalfam\relax#1}}
\font\tenbbb=msbm10 \font\sevenbbb=msbm7 \font\fivebbb=msbm5
\font\ninebbb=msbm9 \font\sixbbb=msbm6 \font\eightbbb=msbm8
\newfam\bbbfam
\def\Bbb#1{{\fam\bbbfam\relax#1}}

\skewchar\ninei='177  \skewchar\eighti='177  \skewchar\sixi='177
\skewchar\ninesy='60  \skewchar\eightsy='60  \skewchar\sixsy='60
\hyphenchar\ninett=-1 \hyphenchar\eighttt=-1 \hyphenchar\tentt=-1


\catcode`@=11
\newskip\ttglue

\def\tenpoint{\def\rm{\fam0\tenrm}
  \textfont0=\tenrm \scriptfont0=\sevenrm \scriptscriptfont0=\fiverm
  \textfont1=\teni  \scriptfont1=\seveni  \scriptscriptfont1=\fivei
  \textfont2=\tensy \scriptfont2=\sevensy \scriptscriptfont2=\fivesy
  \textfont3=\tenex  \scriptfont3=\tenex \scriptscriptfont3=\tenex
  \textfont\itfam=\tenit  \def\it{\fam\itfam\tenit}%
  \textfont\slfam=\tensl  \def\sl{\fam\slfam\tensl}%
  \textfont\ttfam=\tentt  \def\tt{\fam\ttfam\tentt}%
  \textfont\bffam=\tenbf  \scriptfont\bffam=\sevenbf
   \scriptscriptfont\bffam=\fivebf  \def\bf{\fam\bffam\tenbf}%
\textfont\bcalfam=\tenbcal \scriptfont\bcalfam=\sevenbcal
   \scriptscriptfont\bcalfam=\fivebcal  \def\bcal{\fam\bcalfam\tenbcal}%
  \textfont\bbbfam=\tenbbb \scriptfont\bbbfam=\sevenbbb
   \scriptscriptfont\bbbfam=\fivebbb  \def\Bbb{\fam\bbbfam\tenbbb}%
  \tt \ttglue=.5em plus.25em minus.15em
  \normalbaselineskip=12pt
  \setbox\strutbox=\hbox{\vrule height8.5pt depth3.5pt width0pt}%
  \let\sc=\eightrm  \let\big=\tenbig  \normalbaselines\rm}

\def\ninepoint{\def\rm{\fam0\ninerm}
  \textfont0=\ninerm \scriptfont0=\sixrm \scriptscriptfont0=\fiverm
  \textfont1=\ninei  \scriptfont1=\sixi  \scriptscriptfont1=\fivei
  \textfont2=\ninesy \scriptfont2=\sixsy \scriptscriptfont2=\fivesy
  \textfont3=\tenex  \scriptfont3=\tenex \scriptscriptfont3=\tenex
  \textfont\itfam=\nineit  \def\it{\fam\itfam\nineit}%
  \textfont\slfam=\ninesl  \def\sl{\fam\slfam\ninesl}%
  \textfont\ttfam=\ninett  \def\tt{\fam\ttfam\ninett}%
  \textfont\bffam=\ninebf  \scriptfont\bffam=\sixbf
   \scriptscriptfont\bffam=\fivebf  \def\bf{\fam\bffam\ninebf}%
  \textfont\bcalfam=\ninebcal \scriptfont\bcalfam=\sixbcal
   \scriptscriptfont\bcalfam=\fivebcal \def\bcal{\fam\bcalfam\ninebcal}%
  \textfont\bbbfam=\ninebbb \scriptfont\bbbfam=\sixbbb
 \scriptscriptfont\bbbfam=\fivebbb \def\Bbb{\fam\bbbfam\ninebbb}%
  \abovedisplayskip=10.5pt plus 3pt minus 7.5pt
  \belowdisplayskip=10.5pt plus 3pt minus 7.5pt
  \abovedisplayshortskip=0pt plus 2.5pt
  \belowdisplayshortskip=6pt plus 2.5pt minus 3pt
  \tt \ttglue=.5em plus.25em minus.15em
  \normalbaselineskip=11pt
  \setbox\strutbox=\hbox{\vrule height8pt depth3pt width0pt}%
  \let\sc=\sevenrm  \let\big=\ninebig  \normalbaselines\rm}

\def\eightpoint{\def\rm{\fam0\eightrm}
  \textfont0=\eightrm \scriptfont0=\sixrm \scriptscriptfont0=\fiverm
  \textfont1=\eighti  \scriptfont1=\sixi  \scriptscriptfont1=\fivei
  \textfont2=\eightsy \scriptfont2=\sixsy \scriptscriptfont2=\fivesy
  \textfont3=\tenex  \scriptfont3=\tenex \scriptscriptfont3=\tenex
  \textfont\itfam=\eightit  \def\it{\fam\itfam\eightit}%
  \textfont\slfam=\eightsl  \def\sl{\fam\slfam\eightsl}%
  \textfont\ttfam=\eighttt  \def\tt{\fam\ttfam\eighttt}%
  \textfont\bffam=\eightbf  \scriptfont\bffam=\sixbf
   \scriptscriptfont\bffam=\fivebf  \def\bf{\fam\bffam\eightbf}%
  \textfont\bcalfam=\eightbcal  \scriptfont\bcalfam=\sixbcal
   \scriptscriptfont\bcalfam=\fivebcal  \def\bcal{\fam\bcalfam\eightbcal}%
  \textfont\bbbfam=\eightbbb \scriptfont\bbbfam=\sixbbb
   \scriptscriptfont\bbbfam=\fivebbb  \def\Bbb{\fam\bbbfam\eightbbb}%
  \abovedisplayskip=9pt plus 2pt minus 6pt
  \belowdisplayskip=9pt plus 2pt minus 6pt
  \abovedisplayshortskip=0pt plus 2pt
  \belowdisplayshortskip=5pt plus 2pt minus 3pt
  \tt \ttglue=.5em plus.25em minus.15em
  \normalbaselineskip=9pt
  \setbox\strutbox=\hbox{\vrule height7pt depth2pt width0pt}%
  \let\sc=\sixrm  \let\big=\eightbig  \normalbaselines\rm}

\def\tenbig#1{{\hbox{$\left#1\vbox to8.5pt{}\right.\n@space$}}}
\def\ninebig#1{{\hbox{$\textfont0=\tenrm\textfont2=\tensy
  \left#1\vbox to 7.25pt{}\right.\n@space$}}}
\def\eightbig#1{{\hbox{$\textfont0=\ninerm\textfont2=\ninesy
  \left#1\vbox to 6.5pt{}\right.\n@space$}}}

\tenpoint

\input epsf

\def\d{{\rm d}}

\def\E{{\rm E}}
\def\L{{\rm L}}
\def\Linf{{\rm L}^\infty}
\def\P{{\rm P}}
\def\Q{{\rm Q}}
\def\Pth{{\rm P}^\theta}
\def\Qth{{\rm Q}^\theta}

\def\bcalF{{\bcal F}}

\def\card{\mathop{\rm card}\nolimits}

\def\Dtwo#1{D^{(2)}_{#1}}
\def\bbbN{{\Bbb{N}}}
\def\bbbP{{\Bbb{P}}}
\def\bbbQ{{\Bbb{Q}}}
\def\bbbE{{\Bbb{E}}}
\def\bbbR{{\Bbb{R}}}

\def\half{{1\over 2}}

\parindent=4true mm
\def\noas{\noalign{\smallskip}}
\def\Example#1{{\smallbreak\noindent{\bf Example #1}}}
\def\Remark#1{{\smallbreak\noindent{\bf Remark #1}}}

\def\Proof{{\smallbreak\noindent{\sl Proof}}}
\def\Step{\smallbreak\noindent{S{\er TEP} }}

\def\sqr#1#2{{\vcenter{\vbox{\hrule height#2pt
    \hbox{\vrule width#2pt height#1pt \kern#1pt
      \vrule width#2pt}
      \hrule height#2pt}}}}

  \def\endsqr{\mathchoice\sqr{7}{.4}\sqr{7}{.4}\sqr{2.1}{.3}\sqr{2.1}{.3}}
  \def\endproof{\hfill$\endsqr$\medskip}

\def\LscapeMatlabPic#1#2#3{ 
  \goodbreak\topinsert
  \vskip 2.3true cm
  \hskip 1.9true cm
  \hbox{\ \ \epsfxsize=7.6true cm \epsfbox[100 0 400 500]{#1} }
  \vskip -5.2true cm
#3
\vskip -1true mm
  \endinsert }

\def\LscapeSvenPic#1#2#3{ 
  \goodbreak\topinsert
  \vskip 2.3true cm
  \hskip 1.9true cm
  \hbox{\ \ \epsfxsize=7.6true cm \epsfbox[50 0 400 500]{#1} }
#3
\vskip -1true mm
  \endinsert }

\font\twr=cmr12
\font\er=cmr8
\font\sc=cmr9
\pageno=1

\vglue 1true cm

\def\bX{{\bf X}}
\def\bone{{\bf 1}}
\def\bP{{\bf P}}
\def\bPL{{\bf P}^{}_{\!\calL}}
\def\bPTh{{\bf P}^{}_{\!\Theta}}

  \def\calL{{\cal L}}
  
 \def\calX{{\cal X}}

\def\frac#1#2{{#1\over#2}}

\def\Colon{{\colon\,}}
\def\parskip{\medskip}
\def\sectitle#1{\bigskip\centerline{\bf \secnum. #1}\medskip\noindent}

\def\seqno#1{\eqno(\secnum.#1)}
\def\sno#1{{\secnum.#1}}
\def\eqn#1{(\sno#1)}
\def\Step{\smallbreak\noindent{S\er TEP} }
\def\nh{\par\noindent\hangindent=5 true mm\hangafter=1}
\def\snh{\smallskip\nh}

\font\twr=cmr12
\font\er=cmr8
\font\sc=cmr9

\def\etal{{\sl et al.}}

\def\noas{\noalign{\smallskip}}

\def\tJ{\tilde{J}}
\def\tR{\widetilde{R}}

\catcode`@=11
\newskip\ttglue

\baselineskip=15pt
\vglue 1 cm
{\centerline{\twr TWO LILYPOND SYSTEMS OF FINITE LINE-SEGMENTS}

}
$$\matrix{
\hbox{D. J. D{\er ALEY}}\cr
 \noas
 \hbox{\sl Centre for Mathematics and its Applications}\cr
 \hbox{\sl The Australian National University}\cr
 \hbox{\sl Canberra, ACT 0200}\cr }
\quad
\matrix{
\hbox{S{\er VEN} E\er BERT}\cr
 \noas
 \hbox{\sl Institut f\"ur Stochastik}\cr
 \hbox{\sl Karlsruher Institut f\"ur Technologie,}\cr
 \hbox{\sl 76128 Karlsruhe, Germany}\cr}$$
$$\matrix{
\hbox{G{\er\"UNTER} L\er AST}\cr
 \noas
 \hbox{\sl Institut f\"ur Stochastik}\cr
 \hbox{\sl Karlsruher Institut f\"ur Technologie,}\cr
 \hbox{\sl 76128 Karlsruhe, Germany}\cr}$$

\bigskip
\centerline{\bf Abstract}
{\narrower
\medskip
The paper discusses two models for non-overlapping finite
line-segments constructed via the lilypond protocol, operating here
on a given array of points $\bP=\{\P_i\}$ in $\bbbR^2$ with which
are associated directions $\{\theta_i\}$. 
At time 0, for each and every $i$, a line-segment $\L_i$
starts growing at unit rate around the point $\P_i$ in the direction
$\theta_i$, the point $\P_i$ remaining at the centre of $\L_i$; each
line-segment, under Model 1, ceases growth when one of its ends hits another
line, while under Model 2, its growth ceases either when one of its ends
hits another line, or when it is hit by the growing end of some other line.

The paper shows that these procedures are well-defined and gives constructive
algorithms to compute the half-lengths $R_i$ of all $L_i$.
Moreover it specifies assumptions under which stochastic versions, i.e.\ 
models based on point processes, exist.  Afterwards it deals with the
question as to whether there is percolation in Model 1. The paper concludes
with a section containing several conjectures and final remarks.

\vfill

{\eightpoint  \noindent
2000 Mathematics Subject Classification:
 Primary 60D05, Secondary 62M30, 60G55

\parskip
\noindent{\bf Acknowledgement.} 
\baselineskip=12pt
DJD first learnt of a slightly different
line-segment model defined via the lilypond protocol from Marianne M\aa nsson,
when she was at Chalmers University, Goteborg; she noted the existence of
cycles in what is our Model 1.  His subsequent work was done
in part while visiting the Mittag-Leffler Institute, Aarhus University,
\'Ecole Normale Sup\'erieure Paris, and Karlsruhe Institute of Technology
(several times).
SE's and GL's work was supported by the
German Research Foundation (DFG) through the research unit 
"Geometry and Physics of Spatial Random Systems" 
under the grant LA 965/6-1.
We thank Christian Hirsch for telling us of his work before submission.

}}

\vfill\eject


\def\secnum{1}
\sectitle{Introduction and models}%
Suppose given a locally finite set $\bP = \{\P_i\} = \{(x_i,y_i)\}$
of points in the plane; associate with each point a
direction $\theta_i\in[0,\pi)$.  Write $\Pth_i=(\P_i,\theta_i)$ and $\bPTh =
\{\Pth_i: \P_i\in\bP\}$.  When no two directions coincide the doubly-infinite lines
$\L_i^\infty$, $\L_j^\infty$ say, drawn through $\P_i$, $\P_j$
with respective directions $\theta_i$, $\theta_j$
meet in some point $\P_{ij}$ say, so 
$\P_{ij} = g(\Pth_i, \Pth_j)$ for some function $g$.
A {\sl lilypond system of line-segments\/} is
constructed by growing line-segments $\{\L_i\}$, one
through each point $\P_i$ in direction $\theta_i$, their growth starting
at the same time and at the same rate for each
segment, in such a way that $\L_i$ always has $\P_i$ as its
mid-point. We use $\bPL$ to denote the family $\{(\Pth_i, R_i)\}$,
where $R_i$ is the half-length (`Radius') of the
line-segment $\L_i$ (we describe shortly how $R_i$ is determined).

Under {\sl Model 1}, any given line-segment ceases
growth when one of its ends reaches any other
line-segment.  Thus the line-segment $\L_i$ grown through
$\P_i$ stops growing when for the first time it reaches the
point of intersection $\P_{ij}$ for some $j\ne i$ for which $\L_j$
has reached $\P_{ij}$ earlier; if there is no such $j$ then $\L_i$
grows indefinitely. 

Under {\sl Model 2} any given line-segment ceases growth at the first instant
either that one of its ends touches another line-segment or that it is
touched by some other line-segment. 
In contrast to Model 1 an infinite line-segment can exist only if it does
not touch any other line nor does any other line touch it.

A third system of line-segments based on $\bPTh$ leads to the so-called
Gilbert tessellation; its growth resembles Model 1 except that the two parts
of the line, one each side of $\P_i$, each stops its growth independently by
touching another line (Noble (1967) described this construction, basing his
exposition on E.N.\ Gilbert's manuscript `Surface Films of Needle-Shaped 
Crystals').

Models 1 and 2 with their different growth-stopping rules produce rather
different families of line-segments (see e.g.\ Figures 3a and 4): Model 1
produces a `denser' family of line-segments.
To describe some of these differences we use the
ideas of {\sl neighbours}, {\sl clusters}, {\sl doublets} and {\sl cycles}.
Two line-segments are {\sl neighbours\/} when they touch each other.  A
family or set $C$ of
line-segments forms a {\sl cluster\/} when (a) every line-segment in $C$
has a neighbour in $C$, and (b) to every pair of line-segments in $C$,
$\L_0$ and $\L_n$ say, we can find $\{\L_i,\ i=1,\ldots,n-1\} \subseteq C$
such that $\L_{j-1}$ and $\L_j$ are neighbours for $j=1,\ldots,n$.  A cluster
$C$ is finite or infinite according to the number of line-segments it contains.
For Model 2, two line-segments constitute a
{\sl doublet\/} if they are neighbours and of the same size.
Finally, for Model 1, for any given integer
$r=3,4,\ldots,$ the line-segments $\L_1,\ldots,\L_r$ constitute an
$r${\sl-cycle of neighbours\/} (an $r$-cycle for short)
if each of the $r$ pairs $(\L_r, \L_1)$ and $(\L_i, \L_{i+1})$,
$i=1,\ldots,r-1$, consists of neighbours. If we assume all clusters to be finite there exist one--one correspondences between clusters and cycles for
Model 1, and clusters and doublets for Model 2.

General lilypond systems of germ--grain models in $\bbbR^d$, of
points and hyperspheres (we call these {\sl standard lilypond models}),
were introduced in H\"aggstr\"om and Meester (1996) and (with numerical work)
in Daley, Stoyan and Stoyan (1999) (= [DSS]) and Daley, Mallows and Shepp
(2000) (= [DMS]); they have been considered further in
Daley and Last (2005) (= [D\&L]), Heveling and Last (2006), and Last and
Penrose (2012). A space-time version with general convex full-dimensional
grains has recently been developed in Ebert and Last (2013).   
Earlier versions of the model exist in the physics
literature under the name ``touch-and-stop model'' (Andrienko, Brilliantov and
Krapivsky, 1994) where the exact 1-dimensional model and solution of [DMS]
were anticipated; both papers have further distinct material.  In contrast to
those systems, the present paper explores aspects of such a system in which the
`grains' are of lower dimension than the space in which they and the `germs'
are located. Models 1 and 2 both incorporate the idea of being `growth-maximal'
in some way: for Model 1 a grain stops growing so soon as one of its
`growth-points' is impeded; for Model 2 a grain stops growing so soon as it
touches or is touched by any other grain.  Thus, both models can be regarded as
`natural' lower-dimensional analogues of the original point-and-hypersphere
standard lilypond models.  Model 2 can be viewed as the limit as $e\uparrow 1$
of a full dimensional germ--grain model in $\bbbR^2$ with randomly oriented
elliptical grains of eccentricity $e$.  

The paper proceeds as follows.  First we give some basic examples of the
Models to get some feel for the behaviour of the growth process.  Section 3
details an algorithm that constructs Model 1 for finite point sets, with
illustrations of realizations from Poisson distributed germs and uniformly and
independently distributed directions.  This algorithm is the first step towards
understanding the Models in a more formal setting in Sections 4 and 5 where
we discuss their existence and uniqueness based on locally finite point sets:
Section 4 has formal definitions that correspond to our intuitive descriptions.
In Section 5 we establish lilypond models based on a broad class of marked
point processes. Under the additional assumption of stationarity we prove in
Section 6 the absence of percolation in Model 2.
Section 7 contains some discussion and further results. In particular we
provide arguments supporting our view that there is no percolation
in Model 1 (i.e.\ it does not contain an infinite cluster).

\LscapeMatlabPic{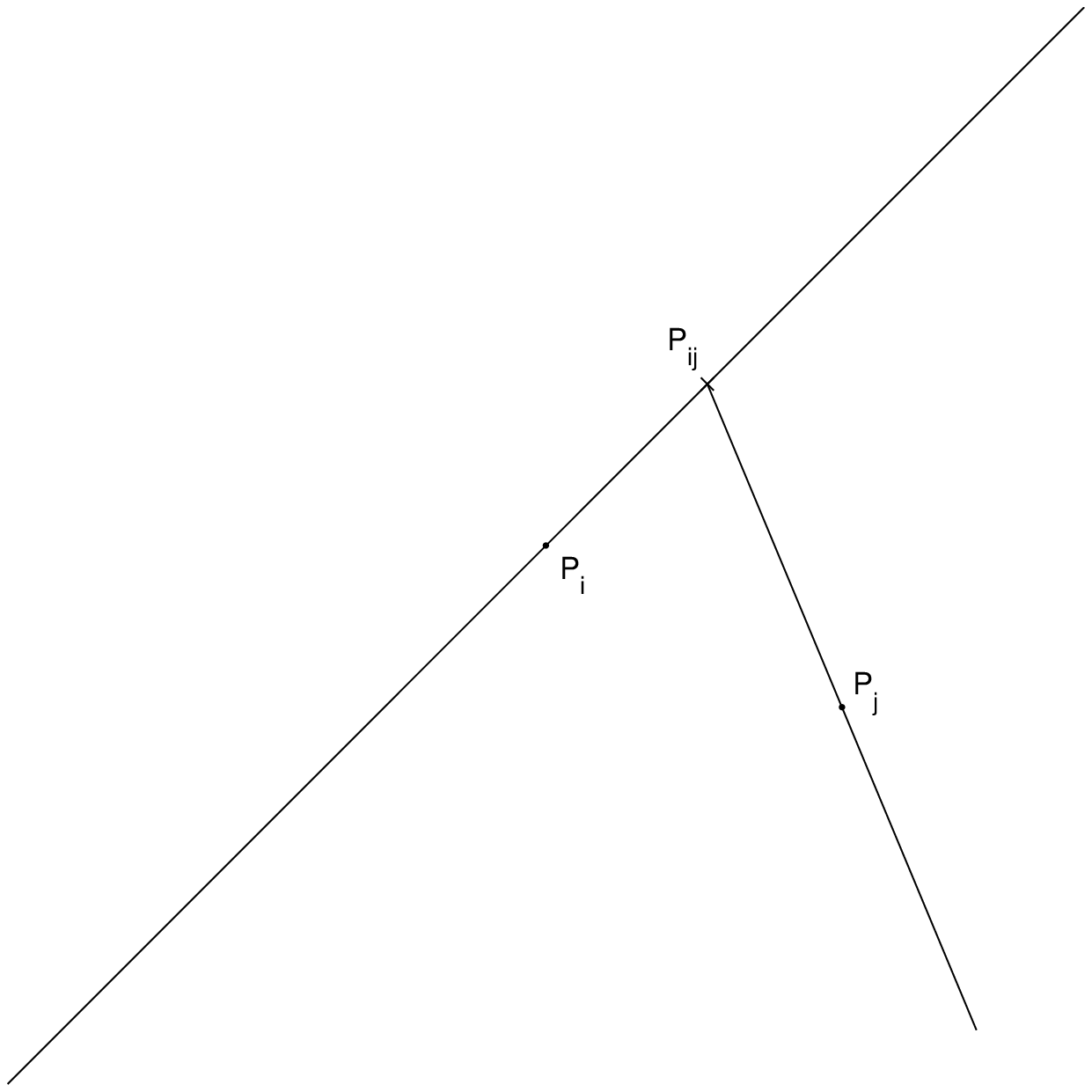}{}{{\narrower\narrower\narrower\narrower
\ninepoint\smallskip
Figure 1. Lilypond line-segments grown through points P$_i$, P$_j$,
meeting in P$_{ij}$. 
\medskip}}

\def\secnum{2}
\sectitle{Basic notation and simple examples}%
Let $d(\P' ,\P'')=|\P'-\P''|$ denote the euclidean distance between two points
$\P' $, $\P''$ in $\bbbR^2$.  We suppose given a set $\bP$ of $n+1$
points $\bP$ and associated directions (in $[0,\pi)$)
 $$\Pth_i = (\P_i, \theta_i) = \big((x_i, y_i), \theta_i\big)
  \qquad(i=0,1,\ldots,n); \seqno{1}$$
let $\bPTh$ denote such a finite family of $\Pth_i$ as in Section 1.
Our analysis mostly uses the distances
 $$d_{ij} := d(\P_i, \P_{ij})\quad\hbox{ and }\quad
d_{ji} := d(\P_j, \P_{ij}),\quad \theta_i\neq\theta_j
\seqno{2}$$
which, for lines growing about centres $P_i$ at unit rate in directions 
$\theta_i$ , represent the times
they need to grow from their germs at $\P_i$ and $\P_j$ to reach
their intersection point $\P_{ij}$.  In the exceptional case that
$\theta_i=\theta_j$, either $P_j$ lies on
the infinite line through $P_i$ with direction $\theta_i$ and we define
$d_{ij}=d_{ji}:=\half d(P_i,P_j)$, i.e.\ the distance between $\P_i$ and the 
midpoint of $P_i$ and $P_j$; else the corresponding lines have an empty
intersection and we set $d_{ij}=d_{ji}:=\infty$.
Then because growth of a line is terminated by touching another line, the
half-segment length $R_i$ must be $D_i^\infty$-valued, where
$D_i^\infty=D_i\cup\{\infty\}$ and
 $$ D_i = \{d_{ij} : d_{ij}>d_{ji}\}. \seqno3$$
We also use  $m_{ij} = \max\{d_{ij}, d_{ji}\} = m_{ji}\,$; these appear in
our discussion of both Models 1 and 2, more notably in the latter because
there the half-segment length $R_i^{(2)}$ is $D_i^{(2),\infty}$-valued, where
now $D_i^{(2),\infty} = D_i^{(2)}\cup\{\infty\}$ and
$D_i^{(2)} := \{m_{ij}: j\ne i\}$.

To obviate the need to refer to exceptional cases assume that all finite
distances $d_{ij}$ are different as in Condition D below
(as a contrary example, using Model 1, if our points were on a lattice
and we restricted growth to lines joining lattice points, Condition D would
be violated frequently and our arguments would be strewn with extra cases).

\proclaim
Condition 2.1 {\rm(Conditions D)}.
A locally finite marked point set $\bPTh$ satisfies\/ {\rm Conditions D} when
all pairwise distances $d_{ij},\ i\ne j$ that are finite, are mutually
disjoint.

Note that in general the occurrence of parallel lines is not excluded by this
condition. As an interesting extreme case we may consider models with only
two different directions.

\Example{1} ({\sl Lilypond line-segment system on two points}).
The simplest nontrivial case consists of two points and their associated
directions, $\bPTh = \{\Pth_i,\,\Pth_j\}$ say. To avoid trivialites we assume
$\theta_i\neq\theta_j$.  When two
line-segments grow in a lilypond system based on such $\bPTh$, the point
$\P_{ij}$ is reached first by the line starting from the point nearer to
$\P_{ij}$, $\P_i$ say, while the line starting from $\P_j$ stops
growing when it reaches $\P_{ij}$ where it touches the line-segment through
$\P_i$ that continues growing indefinitely
(Condition D excludes the possibility that both line-segments are finite and
of the same length).
From \eqn{2}, the finite line-segment
is of half-length $m_{ij}=\max\{d_{ij},d_{ji}\}$.  Specifically,
if $d_{ij} = m_{ij}$, then $R_i = d_{ij}$ finite,
and $R_j=\infty$ (i.e.\ $\L_j=\Linf_j$).

Computationally, the simplest case arises when $\P_0$ is at the origin, 
$\L_0$ is aligned with the $x$-axis, and $\P_1$ is the point of a unit-rate
Poisson process closest to the origin.  The probability density of
$m_{ij}$ is found in Daley \etal\ (2014). \endproof

\Example 2 ({\sl Lilypond line-segment systems on three points}).
Suppose given the set of three marked points $\bPTh = \{\Pth_0,\,
 \Pth_1,\, \Pth_2\}$; apply the lilypond protocol with Model 1.
To exclude exceptional cases assume that no two lines are parallel, i.e.
$\theta_0\neq\theta_1\neq\theta_2\neq\theta_0$.
Because of this, 
a sketch readily shows that some or all of the triangle $\Delta_{012}$ say,
whose vertices are the intersection points $\P_{01}$, $\P_{12}$ and
$\P_{20}$ of the infinite lines $\L_i^\infty$, must also be part of the
line-segments  constructed as a lilypond system, with at most one $\L_i$ of
infinite length.

\LscapeMatlabPic{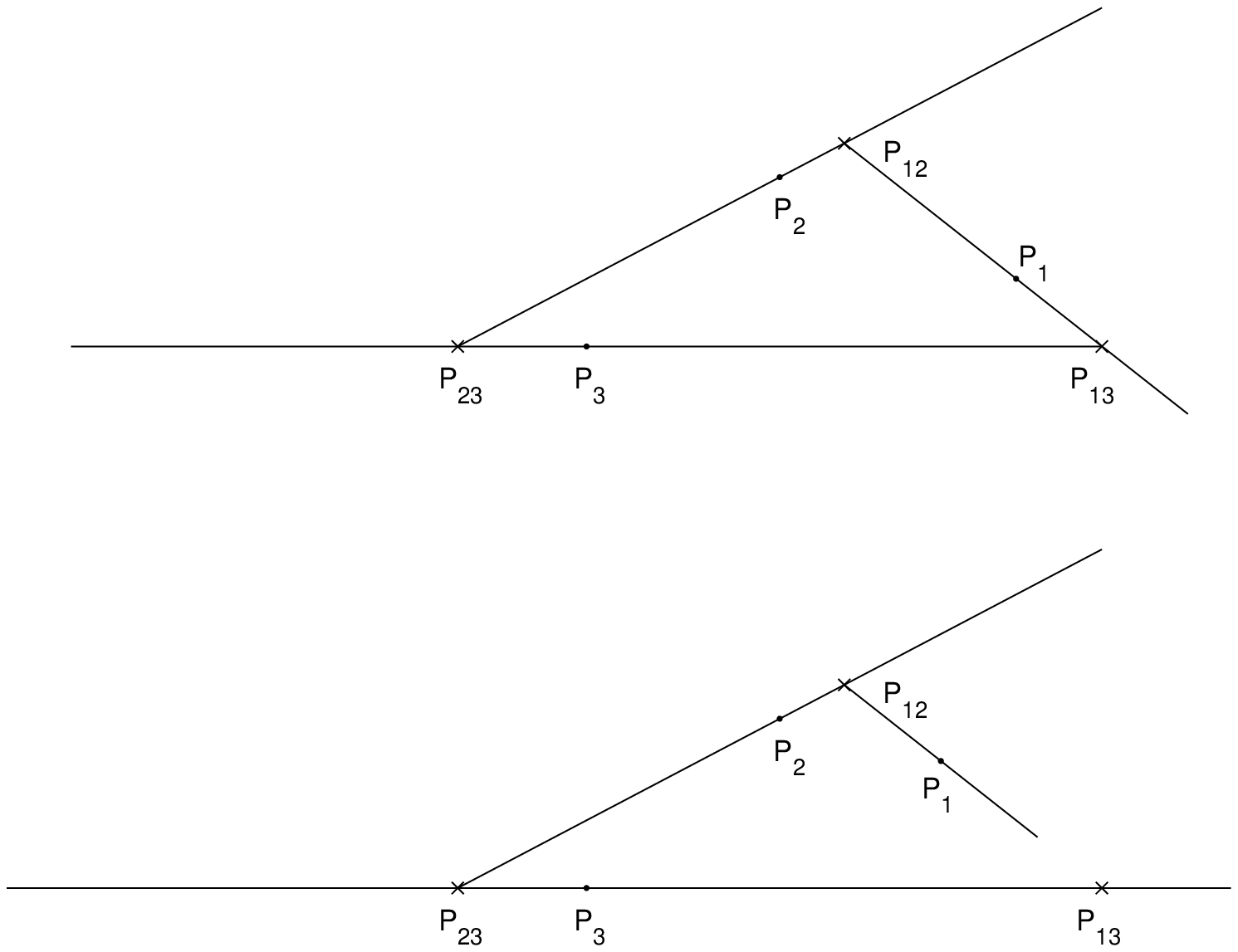}{}{{\narrower\narrower\narrower\narrower
\narrower\ninepoint\smallskip
\noindent
Figure 2.
Two lilypond line-segment models grown through three points: all line-\break
segments finite (left), one infinite (right);
Model 1 (upper), Model 2 (lower).\medskip}}

Recall from around \eqn3 that each half-segment length $R_i$ is
$D_i^\infty$-valued.  For a three-element set $\bPTh$,
each $D_i$ can have at most two elements,
and the union of all three sets must contain exactly three elements.  But for
$R_i$ to be finite, $D_i$ must be non-empty, so for all three $R_i$ to be
finite we cannot have $\{d_{ij} < d_{ji}\hbox{ (all }j\ne i)\}$ for any
$i=0,1,2$.  Defining the sets $A_{ij} = \{d_{ij} < d_{ji}\}$, and recognizing
that (in a space of realizations of 3-element sets $\bPTh$) $A_{ij}\cup A_{ji}$
is the whole space $A$ say, we can write (omitting $\cap$ from
set-intersections in the second and third lines below)
$$\eqalignno{
A &= \bigcap_{0\le i<j\le 2}A_{ij}\cup A_{ji}\,
\,=\, (A_{01}\cup A_{10}) \cap (A_{12} \cup A_{21}) \cap (A_{20}\cup A_{02})\cr
&=
 A_{01}A_{12}A_{20} \cup A_{10} A_{21} A_{02} \cup
A_{01}A_{02}(A_{12} \cup A_{21}) \cup A_{10} A_{12}(A_{20}\cup A_{02})
\cup A_{20}A_{21}(A_{01}\cup A_{10})\cr
  &= A_{01}A_{12}A_{20} \cup A_{10}A_{21}A_{02}
   \cup A_{01}A_{02} \cup A_{10}A_{12} \cup A_{20} A_{21}\,.
&\eqn6\cr}$$
The last three set-intersections in \eqn6 imply
$R_i=\infty$ $(i=0,1,2)$ respectively, while the first two terms of \eqn6
detail two distinct sets of conditions, of which one set necessarily holds if
all three $R_i$ are finite.  Conversely, supposing all $R_i<\infty$, we can
without loss of generality assume $R_0 = \min\{m_{01}, m_{12}, m_{20}\},{}=
d_{01}$ say, implying that $R_1\ge d_{10}$ and, being finite, it must equal
$d_{12}$.  This in turn implies that $R_2\ge d_{21}$ and thus it must equal
$d_{20}$, with $R_0>d_{02}$.  Hence, $A_{10}A_{21}A_{02}$ holds, and
 $\L_0$, $\L_1$ and $\L_2$ form a 3-cycle.   Similarly, still with $R_0 =
\min\{m_{01}, m_{12}, m_{20}\}$ but now${}=d_{10}$, all $R_i$ finite now
implies that $A_{01}A_{12}A_{20}$ must hold, and there is a 3-cycle.

Figure 2 illustrates two possibilities that arise when all three points of
$\bP$ lie on the sides of $\Delta_{012}$; applying Model 1 leads in the upper
case to a 3-cycle and in the lower case to one infinite line-segment.

When Model 2 is based on the three-point set $\bPTh$, we see
that, even with mutually distinct
directions and the centres $\bP$ all lying on the sides of
$\Delta_{012}$, either every line-segment touches another (and all
are of finite length), or one line-segment is of infinite length
(and touches no other). But in no case can we get a 3-cycle as in Model 1.
The analogue for Model 2 of a cycle in Model 1 is a doublet as for
the standard lilypond model in e.g.\ Daley and Last (2005) and
as defined earlier (see above Example 1; in the formal
language of Definition 4.1(c) below,
two points form a {\sl doublet\/} if they are mutual stopping neighbours).
\endproof

Example 2, like Figures 3a and 4, illustrates a major difference between Models
1 and 2: Model 1 leads to cycles coming from at least three points $\Pth_i$,
while Model 2 yields doublets that come from exactly two points.
Despite apparently similar growth rules, the resulting Models are
topologically different.

However, for clusters, the roles of cycles and doublets are similar in
that in Model 1 (resp.\ Model 2) every finite cluster contains exactly one
cycle (resp.\ doublet), and any infinite cluster that may exist
contains at most one cycle (resp.\ doublet).

For Model 1, Examples 1 and 2 differ in that Example 1
always has a line-segment of infinite length but in Example 2 it is
quite possible for all three line-segments to be of finite length.
Inspection of Figures 3a and 3b suggests that for
$\bPTh$ with a larger number $n$ of marked points,
 the occurrence of a line-segment of infinite
length should be increasingly rare as $n$ increases.

\LscapeMatlabPic{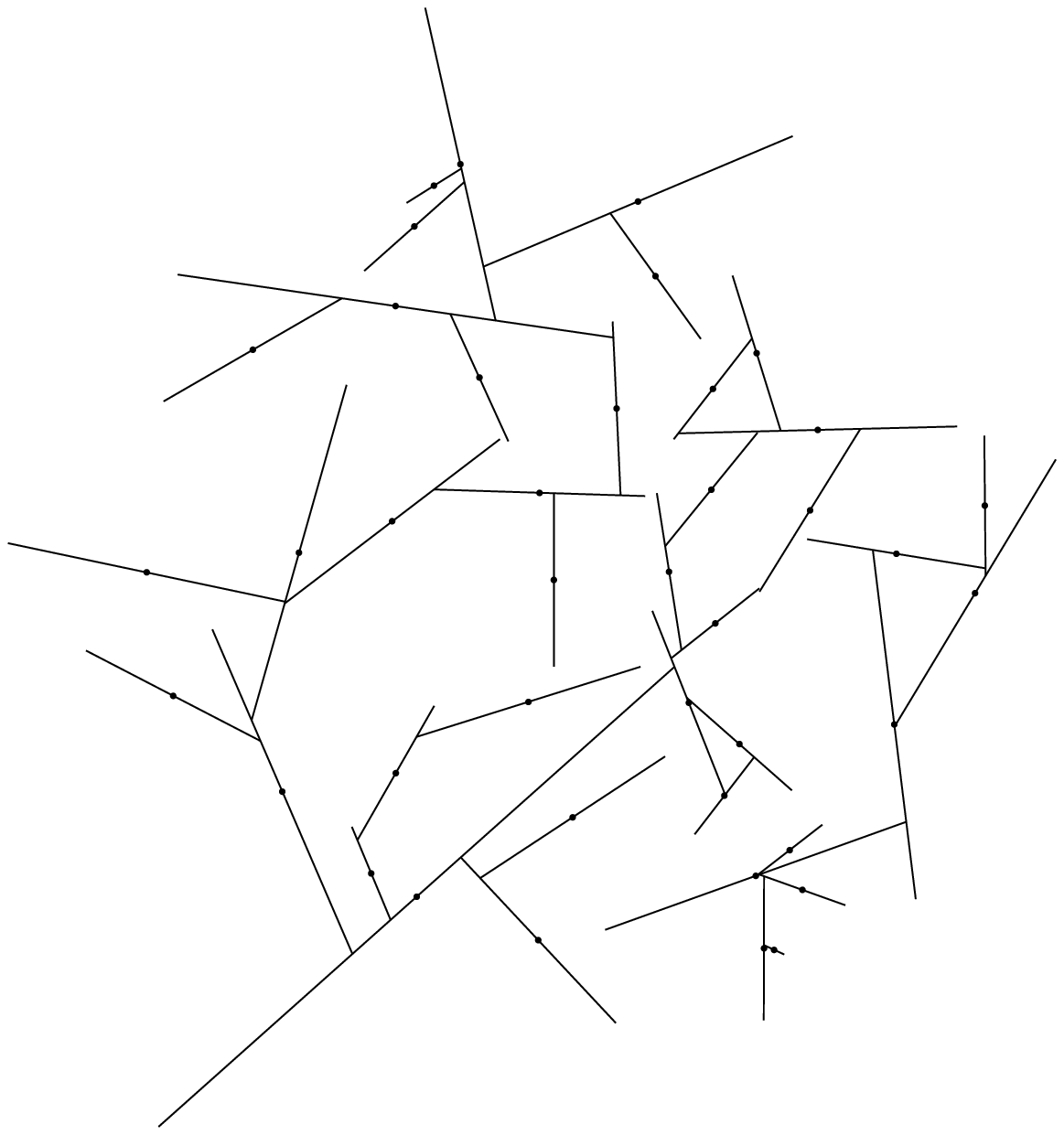}{}{{\narrower\narrower
\ninepoint\medskip
\noindent Figure 3a. All lilypond line-segments grown through 41 Poisson
distributed points $({{\tenpoint\cdot}})$ [Model 1].\bigskip}}

\def\secnum{3}
\sectitle{Solution procedures to find line-segment lengths for finitely many
points}%
We turn to an algorithmic description of Model 1 and briefly sketch
the essentials for Model 2.
The algorithm is generally applicable to a finite marked point
set $\bPTh$.
Given a point $\P_{i_0}$ with index $i_0$, the aim is to identify a
chain of line-segments with mid-points $\P_{i_0},\ldots, \P_{i_{n+r}}$ with
indices $i_0,i_1,\ldots,i_n,\ldots,i_{n+r}$
for which, for $t=0,\ldots,n+r-1$, $\L_{i_t}$ stops growing when it touches
$\L_{i_{t+1}}$ and $\L_{i_{n+r}}$ stops growing when it touches $\L_{i_{n+1}}$
(the chain ends in an $r$-cycle), and $R_{i_t} = d_{i_t, i_{t+1}}$.
The indices are identified sequentially, but we must
allow for the possibility that one $\L_{i_t}$ grows forever; further,
en route from $\P_{i_t}$ while $\P_{i_{t+1}}$ is being found, there may
be branch-chains with indices $j_1,j_2,\ldots\ $.
The strategy underlying the algorithm is similar to that in [DSS]: use a
sequence of lower bounds on $R_i$ to find the earliest time at which the line
$\L_i$ must cease growing.

We now describe an exhaustive algorithm that determines all $R_i$ for a given
finite set $\bPTh$.  What is given below is more efficient and
more informative about the structure of a system of line-segments.

\LscapeMatlabPic{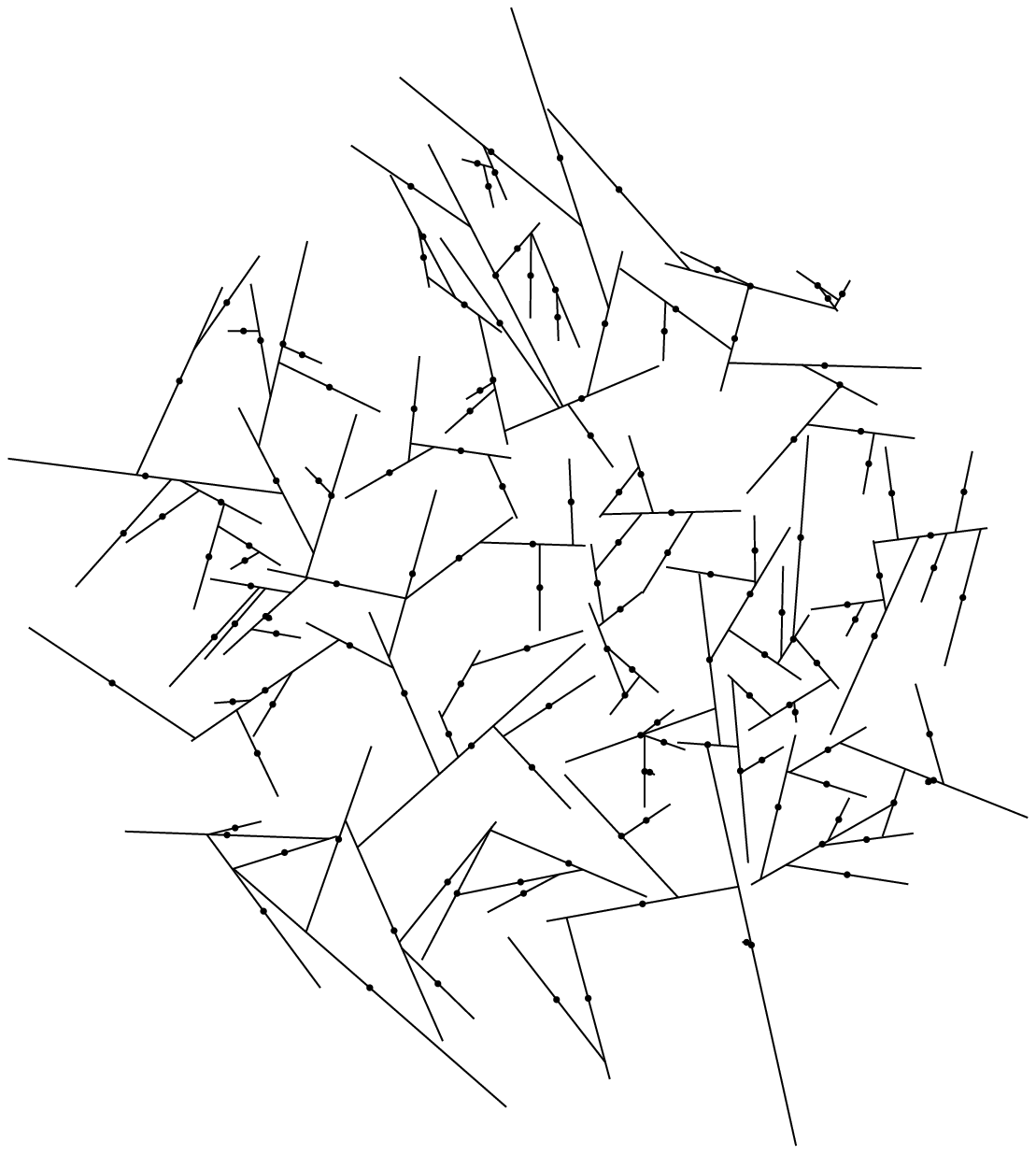}{}{{\narrower\narrower\narrower\narrower
\ninepoint\medskip
\noindent Figure 3b. Same as Figure 3a, but 151 points (innermost 41 points
from Figure 3a).\medskip}}


We have already noted above (2.3) that because
$\L_i$ stops growing by hitting another line-segment $\L_j$ say,
and hence at the intersection-point $\P_{ij}$ as in Example 1,
$R_i$ must be one of the half-lengths in the set $D_i$ defined at (2.3),
implying that $R_i \ge \inf D_i$ provided
$D_i$ is nonempty, else $R_i=\infty$.
If $R_i=d_{ij}$ then as well as $d_{ij}\in D_i$ the line $L_j$ must have
grown at least to $\P_{ij}$, so $R_j>d_{ji}$.  Combining these two facts
implies that
$\{R_i\}$ must satisfy the fixed point relation
 $$R_i = \inf\{d_{ij}: d_{ij}>d_{ji}\hbox{ and } R_j > d_{ji}\}. \seqno1$$
Define $J(i) = \arg\inf\{d_{ij}: d_{ij}>d_{ji}\hbox{ and }R_j>d_{ji}\}$.
Then $R_i = d_{i,J(i)}$, and in terms of the chain $i_0,\ldots,i_{n+r}$
introduced earlier, $J(i_t) = i_{t+1}$ for $t=0,\ldots,n+r-1$ and
$J(i_{n+r}) = i_{n+1}$.  [We digress momentarily to Model 2, for which
 $D_i$ at (2.3) is replaced by the larger set $D_i^{(2)}$ as below
(2.3) and \eqn1 becomes
$$R_i^{(2)} = \inf\{m_{ij} =  \max\{d_{ij}, d_{ji}\} : j\ne i \hbox{ and }
R^{(2)}_j \ge  d_{ji}\}.] \seqno2 $$

Suppose elements $i_0,\ldots,i_t$ of the chain are known; to
identify $J(i_t)=i_{t+1}$ say, we exploit  variants of
\eqn1 and the function $J(\cdot)$.  Write $i=i_t$ and `approximate'
both $R_i$ and $J(i)$ via lower bounds $\tR_j=\inf D_j$ and `trial' elements
$\tJ_q= \arg \inf D_{\tJ_{q-1}}$ for $q=1,2,\ldots,$ with $\tJ_0=i$; strictly,
$\tJ_q = \tJ_q(i)$.  As the
`solution' evolves, the various sets $D_j$ may contract (as potential solutions
$d_{ij}$ are rejected because $R_j < d_{ji}$) and the branch chain
$\tJ_0, \tJ_1,\ldots,$ apart from $\tJ_0 = i$, may also change until
$R_i$ is determined.  The steps below yield both the chain
$i_0,\ldots,i_{n+r}$
and the cycle length $r$.

\LscapeSvenPic{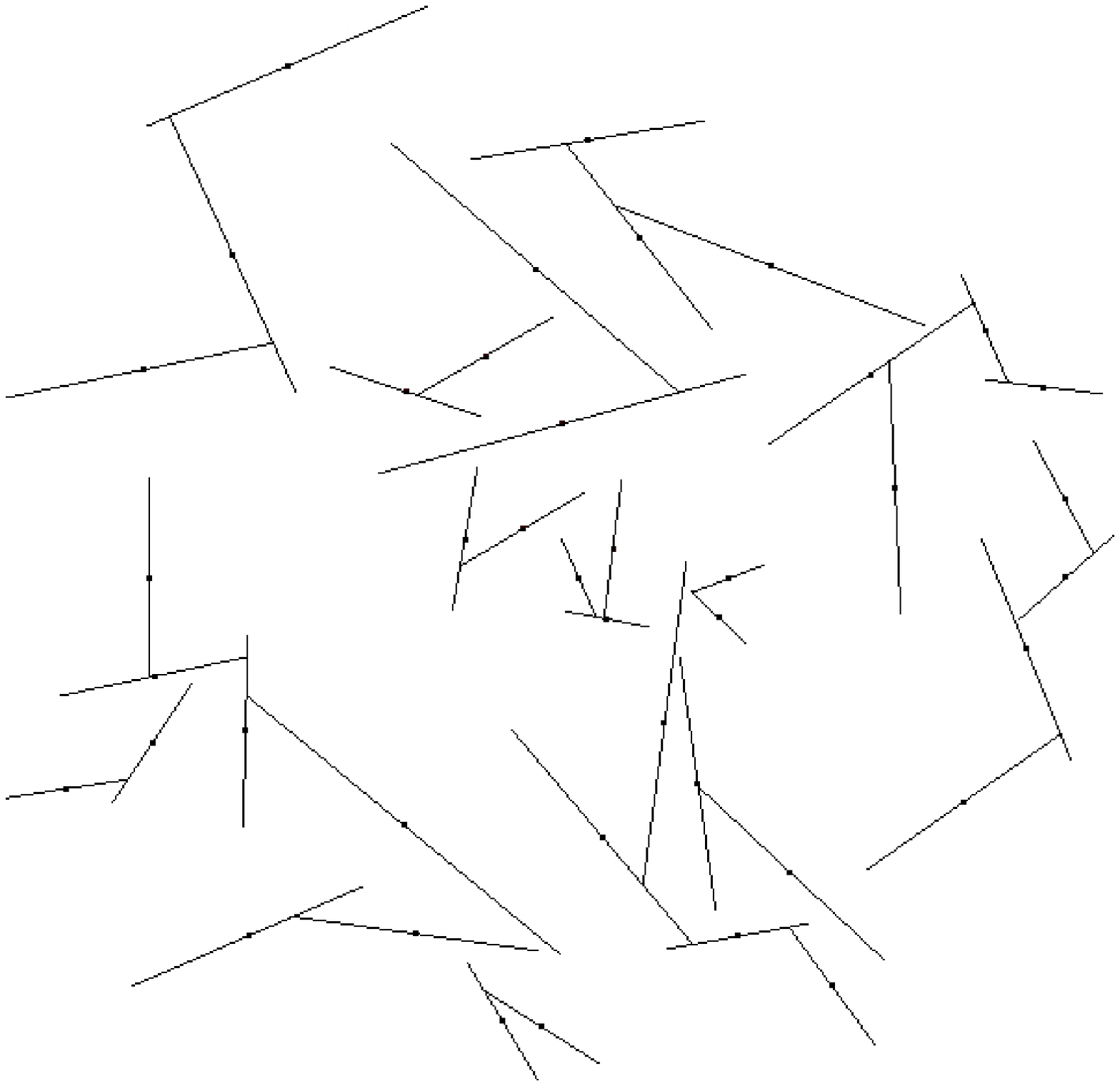}{}{{\narrower\narrower\narrower\narrower
\ninepoint\smallskip
\noindent Figure 4. Model 2 version of Figure 3a.
\medskip}}

\smallskip\noindent{\bf Algorithm 3.1}.  Let the index $i = i_0$ of some point
$\P_{i_0}$ be given; we seek the chain $i_0,i_1,\ldots$ as above, ending
either with an infinite line or an $r$-cycle for some $r$ that is also to be
found.  Set $t=0$.
{\parindent=20pt
\Step 1. Set $q=0$, $\tJ_0=i := i_t$, and construct range-set for $R_i$ viz.\
$D_i = \{d_{ij}: d_{ij} > d_{ji}\}$.
\Step 2. If $D_{\tJ_q}$ is empty, go to 6.4.
Otherwise identify potential stopping index $\tJ_{q+1} := \arg\inf D_{\tJ_q}$
and lower bound $\tR_{\tJ_q} = d_{\tJ_q \tJ_{q+1}}$; set $q\to q+1$.
\item{2.1.} If $q=1$ construct (next) $D_{\tJ_q}$ and repeat Step 2.
\Step 3. If $D_{\tJ_q}$ known go to 3.2; otherwise, construct it.
\item{3.1.} Identify $\tJ_{q+1} := \arg\inf D_{\tJ_q}$, set
$\tR_{\tJ_q} = \inf D_{\tJ_q} = d_{\tJ_q \tJ_{q+1}}$
and go to Step 4.
\item{3.2.} If $R_{\tJ_q}$ known go to Step 5; otherwise go to Step 4.
\Step 4 ({\sl Weak test}\/).  If $\tR_{\tJ_q} < d_{\tJ_q \tJ_{q-1}}$ then $q\to q+1$,
construct $D_{\tJ_q}$ and return to Step 3.1.
\item{4.1.} Otherwise,
$\tR_{\tJ_q} > d_{\tJ_q\tJ_{q-1}}$ so that
$R_{\tJ_{q-1}}$ is found; set $q\to q-1$ and go to Step 6.
\Step 5 ({\sl Strong test}\/).
If $R_{\tJ_q} < d_{\tJ_q \tJ_{q-1}}$ delete
$d_{\tJ_{q-1}\tJ_q}$ from $D_{\tJ_{q-1}}$, $q\to q-1$ and return to
Step 2.
\item{5.1.} Otherwise, $R_{\tJ_q} > d_{\tJ_q\tJ_{q-1}}$ so that
$R_{\tJ_{q-1}}$ is found; set $q\to q-1$ and go to Step 6.
\Step 6. If $q\ge 1$ return to Step 5.
\item{6.1.} Otherwise $R_{i_t} = d_{\tJ_0\tJ_1}$ is found.
If $t=0$ or 1 go to 6.3.
\item{6.2.} If $\tJ_1 = i_{t+1-u}$ for some $u=3,4,\ldots,t$,
then $u=:{}$the cycle length $r$ and Exit. Otherwise,
\item{6.3.} Set $i_{t+1}= \tJ_1 =: J(i_t)$,
$t\to t+1$, and return to Step 1 with new $i=i_t$.
\item{6.4.} $R_{i_t} = \infty$ and no cycle. Exit.
\endproof
} 

\smallskip
\noindent{\bf Algorithm 3.2}.  To find $\{R_i^{(2)}\}$ (i.e.\ Model 2), use the
steps of Algorithm 3.1 with (cf.\ \eqn1 and \eqn2)
 $D_i$ replaced by $D_i^{(2)}$, and $d_{ij}$ by
$\max\{d_{ij}, d_{ji}\}$ as appropriate.
\endproof

We constructed Figures 3a, 3b and 4 using the algorithm described above for
determining all $R_i$ for a given finite set $\bPTh$ in which $\P_0$ is at
the origin, $\L_0$ is aligned with the $x$-axis, $\Pth_1,\ldots,\Pth_n$ are
the $n$ points closest to the origin of a simulated unit-rate marked planar
Poisson process and the directions are i.i.d.\ r.v.s uniform on $(0,\pi)$,
so that Condition D is met a.s.\ (see Section 5).  In this case the
algorithm can be used for the purpose of simulating characteristics of
a family of line-segments under a Palm distribution for $\bPTh$.

\LscapeMatlabPic{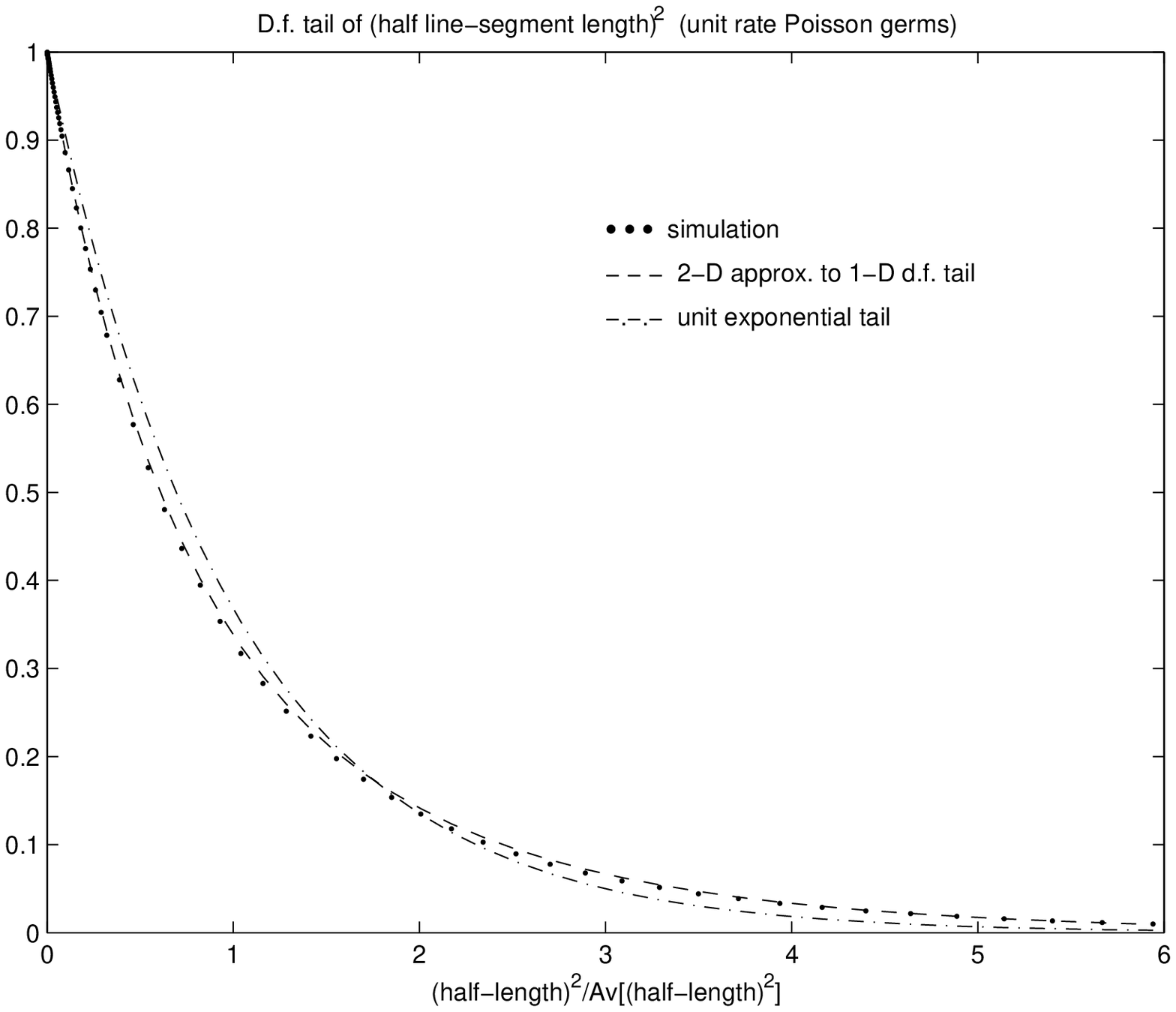}{}{{\narrower\narrower\narrower\narrower
\medskip\ninepoint
\noindent Figure 5.  Tail of the d.f.\ of $[R_1]^2/{\rm Av}[(R_1)^2]$:
observed
\hbox{($\cdot\cdots$)}, transform from exact 1-D tail ($- - -$, [DMS])
exponential with unit mean ($-\cdot - \cdot -$).\medskip}}

We estimated the Palm distribution of a half-line segment $R_i$ in Model 1 by
simulation.  Arguably, it is not $R_i$ but $\pi R_i^2$ that should be used
as a measure of the `space' occupied by a line-segment.  This is borne out by
the closeness of the tail of this distribution to that of the tails of the
`volume' of hyperspheres in the standard lilypond germ--grain models in
 $\bbbR^d$ (see Figure 6 in [DSS] and Figure 5).
The approximate commonality of these distributions is presumably attributable
to the facts that (1) the `germs' $\{\P_i\}$ come from a stationary Poisson
process in the `host' space and (2) the `grains' grow `maximally' as shown by
the fixed-point equations
(here, equations \eqn1 and \eqn2 and, for the radii $r_i$ of hyperspheres in
$\bbbR^d$ in standard models,
 $$r_i = \sup\{x : x + r_j \le d(\P_i, \P_j) \hbox{ (all }j\ne i)\},\seqno3$$
the solution of which satisfies
$d_i := \inf_{j\ne i}\{d(\P_i,\P_j) \} \ge r_i \ge \half d_i$ as in
[DSS]).


\def\secnum{4}
\sectitle{Existence and uniqueness of lilypond line-segment systems}%
To this point we have taken for granted the existence of a
line-segment system generated via the lilypond protocol: when
$\bPTh$ is finite, this follows from Algorithm 3.1.  But when $\bPTh$
is countably infinite, more argument is needed, for which
purpose we exploit the approach in Heveling and Last (2006) (we also take
advantage of the technical Condition D); our
notation builds on what we have already used.

%

The line-segment realization $\bPL=\{(\Pth_i,R_i):\Pth_i\in\bPTh \}$ based on $\bPTh$ satisfies certain properties that can be described in terms of pairs of
lines as in Definition \sno1 below.  To this end, for any $\theta\in[0,\pi)$,
let $u(\theta) = (\cos\theta,\sin\theta)$ denote the unit vector in the
direction $\theta$, so that for any scalar $R\ge 0$, the line-segment of
length $2R$ in direction $\theta$ with mid-point $\P=(x,y)$ is the set
$S(\Pth,R) := \{\P +tRu(\theta):-1\le t\le1\} \,=:\,
[\P-Ru(\theta), P+Ru(\theta)]$; this line-segment has relative interior
$S^0(\Pth,R) := \{\P + t R u(\theta): -1 < t < 1\}$.

\proclaim Definition \sno1.  Let $\bPTh$ be a locally finite marked point
set satisfying Conditions D.  Let $\Pth_i\mapsto R(\bPTh,\,\Pth_i) \equiv
R(\Pth_i)=: R_i$ be any $[0,\infty]$-valued measurable mapping on\/ $\bPTh$
such that for every $\Pth_i\in\bPTh$ the mapping determines line-segments
 $$S_i := S(\Pth_i, R_i) :=
\cases{
\{\P_i + t u(\theta): |t| \le R_i\}
  &if $R_i<\infty$,\cr\noas
\hbox{the line }\{\P_i + tu(\theta_i)\Colon t\in\bbbR\} &if $R_i=\infty$.\cr
 } \seqno{1}$$
When $\theta_i \neq\theta_j$ let $\P_{ij}$ be the point of intersection of
$S(\Pth_i,\infty)$ and $S(\Pth_j,\infty)$, let $d_{ij} = d(\P_i,\P_{ij})$ and
$d_{ji} = d(\P_j,\P_{ij})$.
{\parindent=18pt
\item{\rm(a)} The set $\big\{\big(\Pth_i,\,R_i\big)\Colon
 \Pth_i\in\bPTh\big\}$ is a {\rm hard-segment model} {\rm(HS model)}
(based on $\bPTh$) if for any distinct
$\Pth_i$ and $\Pth_j \in \bPTh$ the line-segments $S_i$ and $S_j$
have disjoint relative interiors.
\item{\rm(b)} Distinct $\Pth_i$ and $\Pth_j \in\bPTh$ in a HS model are
{\rm segment neighbours} if $S_i\cap S_j \neq\emptyset$.
\item{\rm(c)} For segment neighbours $\Pth_i$ and $\Pth_j$ and $k=1,2$,
$\Pth_j$ is a\/ {\rm Type $k$ stopping segment neighbour of $\Pth_i$} when
 $$R_i = \cases{d_{ij}\qquad\qquad \quad\ \ \hbox{if }d_{ij} > d_{ji}\hbox{ and }R_j > d_{ji}
&for $k=1$,\cr\noas
\max\{d_{ij}, d_{ji}\}\quad \hbox{if } R_j \ge d_{ji}
&for $k=2$.\cr}$$
For $k=1,2$, a HS model is {\rm growth-maximal of Type $k$}
(i.e.\ a GMHS model of Type $k$), if every $\Pth_i\in\bPTh$ for which
$R_i<\infty$ has a Type $k$ stopping segment neighbour.
\smallskip
} 

Definition \sno1 is similar to one given in Heveling and Last (2006) for
lilypond systems of the germ--grain models on points and hyperspheres in
$\bbbR^d$; the quantities in (a)--(d) above are direct analogues for
line-segments in the plane but could readily be adapted to systems of
flats in $\bbbR^d$.

The remainder of this section is devoted to establishing the existence and
uniqueness of Models 1 and 2.  We do so by showing that for $k=1,2$, Model $k$
from Sections 2 and 3 is a GMHS model of Type $k$.
Proceeding first via intermediate steps, the major part of the
discussion concerns a given fixed locally finite marked point set
$\bPTh$.  We start with Model 1.

\proclaim Definition \sno2 {\rm(Descending chains)}.
Let $\bPTh$ be a locally finite marked point
set.\hfill\break
{\rm(a)} \ $\bPTh$ has a {\rm descending chain of Type 1}
when it contains an infinite sequence $\{\Pth_0, \Pth_1,\ldots\}$
such that both inequalities in $d_{n-1,n} \ge d_{n,n-1} \ge
d_{n,n+1}$ hold for all $n=1,2,\ldots\,$.\hfill\break
{\rm(b)} \ $\bPTh$ has a {\rm descending chain of Type 2} when it
contains an infinite sequence $\{\Pth_0, \Pth_1,\ldots\}$ such that the
inequality $d_{n,n-1} \ge \max\{d_{n,n+1}, d_{n+1,n}\}$ holds for
all $n=1,2,\ldots\,$.

Here then is the result for Model 1;
notice that the right-hand side of \eqn3 is a generalization of
the right-hand side of (3.1), and that the fixed-point equation $f=T_1f$
is an extension of (3.1).

\proclaim Theorem \sno3.
Let $\bPTh=\{\Pth_i\Colon i=1,2,\ldots\}$ be a locally finite marked point set
satisfying Conditions D and such that $\bPTh$ admits no descending chain
of Type 1. Then there exists a unique GMHS model of Type 1 based on $\bPTh$,
and it is the unique solution for $f\in\bcalF$ of\/ $T_1 f = f$, where
$\bcalF$ is the space of measurable functions $f:\bPTh \mapsto [0,\infty]$,
the operator $T_1:\bcalF\mapsto\bcalF$ is defined by
 $$\eqalignno{
T_1 f(\Pth_i) :&= \inf D_i(f,\bPTh) &\eqn2 \cr
\noalign{\hbox{and}}
 D_i(f,\bPTh) :&= \{d_{ij}: \Pth_j \in \bPTh\setminus\{\Pth_i\},\,
d_{ij} > d_{ji} \hbox{ and } f(\Pth_j)> d_{ji}\}. &\eqn3 \cr }$$


Theorem 4.3 is a consequence of several results given below where we omit the
phrase `of Type 1' (since we deal only with Model 1 until Theorem 4.12), and
we assume that $\bPTh$ satisfies Conditions D and that there is no
descending chain (of Type 1).

Start by noting that a HS function is an element of $\bcalF$
satisfying the requirements of Definition \sno{1(a)}, and a GMHS function is
a HS function satisfying the case $k=1$ of Definition \sno{1(c)}.
Proposition \sno11 below identifies the GMHS function as the unique
fixed point of the operator $T_1: \bcalF\mapsto \bcalF$ defined at \eqn2,
and as usual, in \eqn3, $\inf\emptyset = \infty$.
Immediately, for $f,g\in\bcalF$, if $f\le g$ then $D_i(g,\bPTh) \supseteq
D_i(f, \bPTh)$.  Appeal to \eqn2 proves the following monotonicity
property.

\proclaim Lemma \sno4.
Let $f,g\in\bcalF$ satisfy $f\le g$.  Then $T_1f \ge T_1 g$.

The next property gives a simple condition under which $D_i(f,\bPTh)$ is a
finite set so that the infimum at \eqn2 is attained.

\proclaim Lemma \sno5.
Let $f\in\bcalF$ and $\Pth_i\in\bPTh$ satisfy $T_1f(\Pth_i) < \infty$.
Then there exists $\Pth_j\in\bPTh\setminus\{\Pth_i\}$ such that
$f(\Pth_i) = d_{ij} > d_{ji}$ and $f(\Pth_j) > d_{ji}$.

\Proof. Because $\inf\emptyset = \infty > T_1f(\Pth_i)$,
$D_i(f,\bPTh)$ is a nonempty set.  To show that it is a finite set,
observe that for any nonempty triangle $\P_i \P_{ij} \P_j$,
 $2 m_{ij} = 2\max\{d_{ij}, d_{ji}\} \ge d_{ij} + d_{ji} \ge d(\P_i, \P_j)$
so for any $c> 0$,
 $$\big\{\Pth_j \in \bPTh: c \ge d_{ij} > d_{ji}\big\}
\,\subseteq \big\{\Pth_j\in\bPTh : 2c \ge d(\P_i,\P_j)\big\}; \seqno{4}$$
this last set is finite because $\bPTh$ is locally finite.
Take $c>T_1f(\Pth_i)$.  Then
  $\inf D_i(f,\bPTh) =
  \inf\big\{D_i(f,\bPTh) \cap \{j\Colon c \ge d_{ij} > d_{ji}\}\big\} $.
 But by \eqn4 this last set is finite, so
$\card\big(D_i(f,\bPTh)\big) < \infty$,
and the infimum at \eqn2 must be attained at an element of the set.
\endproof

\proclaim Lemma \sno6.  Let $f\in\bcalF$.  Then $f$ is a HS function
if and only if $f\le T_1f$.

\Proof.   Assume that $f$ is a HS function, and take $\Pth_i\in\bPTh$.
To show that $f(\Pth_i) \le T_1f(\Pth_i)$, we argue by contradiction: 
assume that for some $\Pth_i$, $T_1f(\Pth_i) < f(\Pth_i)$.
This implies first that $T_1 f(\Pth_i) < \infty$, and then by Lemma \sno5
that for some $j$ we have $T_1 f(\Pth_i) = d_{ij}$ and so
 $$f(\Pth_i) > T_1f(\Pth_i) = d_{ij} > d_{ji}
  \quad\hbox{\rm and} \quad f(\Pth_j) > d_{ji}\,. \seqno{5}$$
Then $\P_{ij}$ is interior to both line-segments $S\big(\Pth_i,f(\Pth_i)\big)$
and $S\big((\Pth_j,f(\Pth_j)\big)$, contradicting the HS property at
Definition \sno1 for $f$.

Conversely, assume that $f \le T_1f$, and take $\Pth_i\in\bPTh$ and
$\Pth_j \in\bPTh \setminus \{\Pth_i\}$; we must show that the
relative interiors $S_i^0 := S^0\big(\Pth_i,f(\Pth_i)\big)$
and $S_j^0 := S^0\big((\Pth_j,f(\Pth_j)\big)$ have a void intersection.
If these two line-segments are not parallel, any non-void intersection
$S(\Pth_i, \cdot) \cap S(\Pth_j,\cdot)$ consists of the point $\P_{ij}$ which,
being at distances $d_{ij}$ and $d_{ji}$ from $\Pth_i$ and $\Pth_j$, is not in
$S_i^0 \cap S_j^0$ when
$f(\Pth_i) \le T_1 f(\Pth_i) = d_{ij}$ for which $f(\Pth_j) > d_{ji}$ by
definition of $T_1 f$.  If the two line-segments are parallel, then either
the infinite lines that contain them have no finite point of intersection and
$S_i^0\cap S_j^0 = \emptyset$, or they both lie within the same line, in which
case $d_{ij} = d_{ji}$ which is impossible when Condition D holds.
Thus, $f$ is an HS function.
\endproof

\proclaim Lemma \sno7.
Let $f\in\bcalF$.  Then $f$ is a GMHS function if and only if $T_1f = f$.

\Proof. When $f$ is a GMHS function it is an HS function so it is
enough to show that an HS function for which $f = T_1 f$ is a GMHS
function. Take $\Pth_i \in \bPTh$.  Either $f(\Pth_i) = \infty$ and
$f(\Pth_i) = T_1 f(\Pth_i)$, or $f(\Pth_i)<\infty$.  In this
case, as in the proof of Lemma \sno6, any non-void intersection of
line-segments determined by $\Pth_i$ and $\Pth_j$ consists of the
singleton set $\{\P_{ij}\}$, and such line-segments can have void
intersection of their relative interiors only if $\P_{ij}$ is at an
extremity of one of the segments, so for some $j$ we have $f(\Pth_i)
= d_{ij} = T_1 f(\Pth_j) > d_{ji}$ and $f(\Pth_j) > d_{ji}$.
\endproof

Lemmas 4.5--7 imply that {\sl when a locally finite marked point set
$\bPTh$ satisfies Conditions D, Model 1 generates a family of line-segments}.
It remains to show that such a family is unique.

For use below we note the following corollary as a separate result.

\proclaim Lemma \sno8.  Let $f$ be a GMHS function.  Then $f(\Pth_i) \in
\big\{d_{ij} : d_{ij} > d_{ji}\big\}$ whenever $f(\Pth_i) < \infty$.

Define now a sequence of functions $f_n\in\bcalF$ recursively via
 $$f_0 := 0, \qquad f_{n+1} = T_1 f_n \quad(n=0,1,\ldots),\seqno{6}$$
so that $f_1 = \infty$.
Using Lemma \sno5, $f_0\le f_1$ implies $f_1 \ge f_2 \le f_3 \ge f_4 \le
\cdots,$ while $f_0\le f_2$ and $f_1 \ge f_3$ imply that $f_{2n} \le f_{2n+2}$
and $f_{2n+1} \ge f_{2n+3}$ for all $n\ge 0$.  Then the monotone limits
 $$f := \lim_{n\to\infty} f_{2n}, \qquad g := \lim_{n\to\infty} f_{2n+1}
\seqno{7}$$
are well-defined, and
 $$
\!\!\!\!\!
f_{2n} \le f_{2n+2} \le f \le g \le f_{2n+3} \le f_{2n+1}
 \qquad (n\ge 0).
\!\!\!\!\!
\seqno{8}$$

Our aim now is to show that $f=g$, because \eqn7 and \eqn8 then imply that
$f$ is the unique GMHS function.  First we derive some auxiliary results.

\proclaim Lemma \sno9.  Let $\Pth_i\in\bPTh$ satisfy $f(\Pth_i)
 < \infty$.  Then $f_{2n}(\Pth_i) = f(\Pth_i)$ for all sufficiently large $n$.
Similarly, if $g(\Pth_j) < \infty$ then $f_{2n+1}(\Pth_j) = g(\Pth_j)$ for all
sufficiently large $n$.

\Proof.  The assertions follow from Lemma \sno8 and the fact that in \eqn3
the right-hand side, and hence also the left-hand side, is a finite set.
\endproof

\proclaim Lemma \sno{10}. $T_1f = g$ and $T_1g = f$.

\Proof.  From $f_{2n} \le f$ and Lemma \sno4 it follows that
$f_{2n+1} \ge T_1f$ and hence that $g\ge T_1f$.  Consider
$\Pth_i\in\bPTh$: we want to show that $g(\Pth_i) \le T_1
f(\Pth_i)$. When $T_1f(\Pth_i) = \infty$ it follows that $g(\Pth_i)
= T_1f(\Pth_i)$, so we can assume that $T_1f(\Pth_i) < \infty$. By
Lemma \sno5 there exists $\Pth_j\in \bPTh\setminus \{\Pth_i\}$ such
that $T_1f(\Pth_i) = d_{ij} \ge d_{ji}$ and $f(\Pth_j) > d_{ji}$.
Assume that $f(\Pth_j)=\infty$. Then $f_{2n}(\Pth_j) > d_{ji}$ for
all sufficiently large $n$, and thus
 $$f_{2n+1}(\Pth_i) = T_1 f_{2n}(\Pth_i) \le d_{ij} = T_1f(\Pth_i)$$
for all sufficiently large $n$, implying that $g(\Pth_i) \le T_1f(\Pth_i)$.
Assuming $f(\Pth_j) < \infty$, Lemma \sno9 implies that
$f_{2n}(\Pth_j) = f(\Pth_j) > d_{ji}$ for all sufficiently large
$n$.  This again implies that $g(\Pth_i) \le T_1f(\Pth_i)$.

To show that $T_1g = f$, start from $f_{2n+1} \ge g$ and Lemma \sno8
to deduce that $f_{2n+2} \le T_1g$ and hence $f\le T_1g$.  To show
that $f\ge T_1g$, take $\Pth_i\in\bPTh$ and assume on the contrary
that $f(\Pth_i) < T_1g(\Pth_i)$. Then $f_{2n}(\Pth_i) = f(\Pth_i)$
for all sufficiently large $n$. By \eqn5 there must be
$\Pth_j\in\bPTh\setminus \{\Pth_i\}$ such that
 $$f(\Pth_i) = f_{2n}(\Pth_i) = d_{ij} \ge d_{ji} \quad\hbox{ and }
\quad f_{2n-1}(\Pth_i) > d_{ji}$$ for infinitely many $n$. But then
$g(\Pth_j) > d_{ji}$, implying that $T_1g(\Pth_i) \ge d_{ij} =
f(\Pth_i)$, which contradicts our assumption that $f(\Pth_i) <
T_1g(\Pth_i)$. \endproof

\proclaim Proposition \sno{11}.
The function $f$ is a GMHS function based on $\bPTh$
if and only if\/ $f=g$, in which case $f$ is the unique such GMHS function.

\Proof.  Suppose $f=g$. From Lemma \sno{10}, $T_1f = T_1g = f$, which
implies by Lemma \sno7 that $f$ is a GMHS function.  For any GMHS
function $h$ we must have $T_1h=h$. But $f_0\le h$ by definition of
$f_0$, so $f_{2n} \le h$ for every $n$, and therefore $f \le h$. But
by Lemma \sno4 we then have $T_1f \ge T_1h = h$, and $f=T_1f$ so $f\ge
h$, hence $f=h$.

Conversely, if $f=T_1f$ then Lemma \sno{10} implies that $f=g$.\endproof

Theorem \sno3 is now a consequence of the last proposition and the next.

\proclaim Proposition \sno{12}.
Under the assumptions of Theorem \sno3, $f=g$.

\Proof.
We use the inequality $f\le g$ and Lemma \sno{10} without further reference.
Assume that $\Pth_0\in\bPTh$ satisfies $f(\Pth_0) < g(\Pth_0)$, and let
$\Pth_1\in\bPTh\setminus\{\Pth_0\}$ be such that
 $$f(\Pth_0) = T_1 g(\Pth_0) = d_{01} \ge d_{10}
  \quad\hbox{ and }\quad  g(\Pth_1) > d_{10}\,.$$
Then
 $f(\Pth_1) \le d_{01}$
because otherwise, $g(\Pth_0) = T_1f(\Pth_0) \le d_{01} =
f(\Pth_0)$.

 We also have $f(\Pth_1) < g(\Pth_1)$ because otherwise
we should have
 $f(\Pth_1) = g(\Pth_1)$, so that again $g(\Pth_0) = T_1f(\Pth_0)\le d_{01}
 = f(\Pth_0)$.

Hence, we can repeat all steps to deduce the existence of some $\Pth_2\in\bPTh
\setminus \{\Pth_1\}$ such that
 $$f(\Pth_1) = T_1g(\Pth_1) = d_{12} \ge d_{21} \quad\hbox{ and }
\quad g(\Pth_2) > d_{21},$$
and $f(\Pth_2) \le d_{21}$ and $f(\Pth_2) < g(\Pth_2)$.  Combining
these inequalities yields the relations
 $$d_{01} > d_{10} \ge d_{12} > d_{21},$$
in which the strict inequalities come from the first assumption of Theorem
\sno3.  In particular, $\Pth_2\neq \Pth_1$.  By induction we can construct a
whole sequence $\Pth_0, \Pth_1, \ldots$ of points from $\bPTh$ satisfying
 $$d_{01} > d_{10} \ge d_{12} > d_{21}
  \ge \cdots \ge d_{n-1,n} > d_{n,n-1} \ge \cdots$$
and $f(\Pth_n) = d_{n,n+1}$ for all $n\ge 0$.  In particular then,
$f(\Pth_n) > f(\Pth_{n+1})$, showing that the points $\P_n$ are all
different.  But this means that we have constructed a descending chain of
$\bPTh$ contrary to what is assumed in Theorem \sno3.  Hence there can be no
$\Pth_0 \in\bPTh$ such that $f(\Pth_0) < g(\Pth_0)$. \endproof

We turn now to discuss the existence of Model 2 along the lines of the proof
for Model 1: it is to be understood that the analysis for the remainder of
this section concerns Model 2, and that we should refer to GMHS model and
stopping segment neighbours of Type 2.

\proclaim Theorem \sno{13}.
Let $\bPTh=\{\Pth_i\Colon i=1,2,\ldots\}$ be a locally finite marked point set
satisfying Conditions D and such that $\bPTh$ admits no descending chain of
Type 2.
Then there exists a unique GMHS model of Type 2 based on $\bPTh$, and it is
the unique solution for $f\in\bcalF$ as in Theorem 4.3 of $T_2 f = f$, where
the operator $T_2:\bcalF\mapsto \bcalF$ is defined by
 $$\eqalignno{
 T_2f(\Pth_i) &:= \inf \Dtwo{i}(f,\bPTh), &\eqn{10}\cr
 \Dtwo{i}(f,\bPTh) &:= \big\{\max\{d_{ij},d_{ji}\}\Colon \Pth_j \in
\bPTh\setminus\{\Pth_i\} \hbox{ and } f(\Pth_j) \ge d_{ji}\big\}. &\eqn{11}
\cr }$$

Theorem 4.13 is proved via several intermediate results as for Theorem 4.3,
assuming now that $\bPTh$ satisfies Conditions D and
has no descending chain (of Type 2).

We start with a monotonicity result, proved as for
Lemma \sno4, and the attainment of an infimum, proved as for Lemma \sno5
with $\big\{\Pth_j \in \bPTh: c \ge d_{ij} \ge d_{ji}\big\}$  replaced by
 $\big\{\Pth_j \in \bPTh: c \ge \max\{ d_{ij} , d_{ji}\}\big\}$.

\proclaim Lemma \sno14.
Let $f,g\in\bcalF$ satisfy $f\le g$.  Then $T_2f \ge T_2g$.

\proclaim Lemma \sno15. Let $f\in\bcalF$ and $\Pth_i\in\bPTh$
satisfy $T_2f(\Pth_i) < \infty$. Then there exists
$\Pth_j\in\bPTh\setminus\{\Pth_i\}$ such that $T_2 f(\Pth_i) =
\max\{d_{ij}, d_{ji}\}$ and $f(\Pth_j) \ge d_{ji}$.

The next step we prove via four intermediate results.

\proclaim
Proposition \sno{16}.
Let $f\in\bcalF$.  Then $f$ is a GMHS function if and only if\/
$f=T_2 f$.

\proclaim Lemma \sno17.  Let $f\in\bcalF$ and assume $f = T_2f$.
Then $f$ is a HS function.

\Proof. Suppose that $f(\Pth_j) > d_{ji}$ and $f(\Pth_i) > d_{ij}$
for some $i\neq j$. If $d_{ij}\ge d_{ji}$, then we have $f(\Pth_i) =
T_2 f(\Pth_i) \le \max\{d_{ij}, d_{ji}\} = d_{ij}$ which is a
contradiction. If $d_{ij} < d_{ji}$, we get $f(\Pth_j) = T_2
f(\Pth_j) \le \max\{d_{ij}, d_{ji}\} = d_{ji}$, which again gives a
contradiction.
\endproof

\proclaim Lemma \sno18. Let $f\in\bcalF$ and assume $T_2f = f$. Then
$f$ is a GMHS function.

\Proof. Because of Lemma \sno17 we can assume that $f$ is a HS function.
Take $i\in \bbbN$.  By Lemma \sno15 there exists $j\neq i$ such that
$T_2 f(\Pth_i) = \max\{d_{ij}, d_{ji}\}$ and $f(\Pth_j) \ge d_{ji}$.
We claim that $\Pth_j$ is a stopping neighbour of $\Pth_i$.
We do this by considering four cases:
{\parindent=17pt
\item{(1)} Suppose $f(\Pth_i)=d_{ij}$ and $f(\Pth_j)=d_{ji}$. Since all $d_{ij}$
are different we get
$ f(\Pth_i)=d_{ij}=T_2 f(\Pth_i)=\max\{d_{ij}, d_{ji}\}>d_{ji}=f(\Pth_j)$.
Then by definition $\Pth_j$ is a stopping neighbour of $\Pth_i$.
\item{(2)}
Suppose $f(\Pth_i)=d_{ji}$ and $f(\Pth_j)=d_{ji}$. So
$f(\Pth_i)=d_{ji}=f(\Pth_j)$ holds and the claim follows.
\item{(3)}
Suppose $f(\Pth_i)=d_{ji}$ and $f(\Pth_j)>d_{ji}$. This gives
$f(\Pth_i)=d_{ji}> d_{ij}$ and $f(\Pth_j)>d_{ji}$. Since $f$ is a
HS function this is a contradiction.
\item{(4)}
Suppose $f(\Pth_i)=d_{ij}$ and $f(\Pth_j)>d_{ji}$. By assumption
$T_2f(\Pth_i)= \max\{d_{ij}, d_{ji}\}$.  Since
$f(\Pth_i)=d_{ij}$ we get $T_2f(\Pth_j)\le \max\{d_{ij},d_{ji}\}$.
This yields $f(\Pth_i)\ge f(\Pth_j)$ and the claim follows.
\endproof

} 
\proclaim Lemma \sno19. Let $f\in\bcalF$ and assume $f$ is a
GMHS function. Then $f\ge T_2 f$.

\Proof.
If $f(\Pth_i) < \infty$ then there exists $\Pth_j\in\bPTh\setminus\{\Pth_i\}$
such that $f(\Pth_i) = \max\{d_{ij}, d_{ji}\} $, and $f(\Pth_j) \ge d_{ji}$.
This implies that $T_2 f(\Pth_i) \le \max\{d_{ij}, d_{ji}\} = f(\Pth_i)$.
If $f(\Pth_i)=\infty$ the proposition is satisfied, since
$T_2f(\Pth_i)$ takes values in $[0,\infty]\cup\{\infty\}$.
\endproof

\proclaim Lemma \sno20.
Let $f\in\bcalF$ and assume $f$ is a GMHS function.  Then $f\le T_2f$.

\Proof. Let $i\ge 1$. To show that $f(\Pth_i)\le T_2f(\Pth_i)$ it
clearly suffices to assume that $T_2 f(\Pth_i) <\infty$. By Lemma
\sno14 there exists $j\neq i$ such that $T_2 f(\Pth_i) =
\max\{d_{ij}, d_{ji}\}$ and $f(\Pth_j) \ge d_{ji}$. We examine two
cases, supposing first that
$f(\Pth_j)>d_{ji}$ and $T_2f < f$.  Then
$$ d_{ij}\le \max\{d_{ij}, d_{ji}\}=T_2 f(\Pth_i) < f(\Pth_i)$$
and $f(\Pth_j)> d_{ji}$ which would contradict the fact that $f$ is
a HS function.  Suppose on the other hand that
$f(\Pth_j)=d_{ji}$. Since $f$ is a GMHS function, $\Pth_j$ has
a stopping neighbour $\Pth_k$. In particular
$f(\Pth_j)=\max\{d_{jk},d_{kj}\}$ holds. Since all $d_{lm},\ l\neq m$
are different we must have $i=k$. Therefore the point $\Pth_i$
must be a stopping neighbour of $\Pth_j$. If we assume $f(\Pth_i) >
T_2 f(\Pth_i)$ then
 $$ f(\Pth_i) > T_2 f(\Pth_i) = \max\{d_{ij}, d_{ji}\}\ge
d_{ji}=f(\Pth_j).$$
This would be a contradiction since $\Pth_i$ stops $\Pth_j$.
\endproof

Now define limit functions $f$ and $g$ as for Model 1 at \eqn6 and \eqn7
except that $T_2$ replaces $T_1$.  Then Lemma 4.21 is an analogue of
Lemma 4.9, and the proof of Lemma \sno{22} is as for Lemma \sno{10} except
that $d_{ij}$ is replaced by $\max\{d_{ij}, d_{ji}\}$.

\proclaim Lemma \sno21.
Let $\Pth_i\in\bPTh$ satisfy $f(\Pth_i) < \infty$.
Then $f_{2n}(\Pth_i) = f(\Pth_i)$ for all sufficiently large $n$.
Similarly, if $g(\Pth_j) < \infty$ then $f_{2n+1}(\Pth_j) = g(\Pth_j)$
for all sufficiently large $n$.

\proclaim Lemma \sno22. $T_2 f = g$ and $T_2 g = f$.

To prove the next proposition mimic the proof of
Proposition \sno11.

\proclaim Proposition \sno{23}.  The function $f$ is a GMHS function
based on $\bPTh$ if and only if $f=g$, in which case $f$ is the
unique such GMHS function.

Theorem \sno13 is now a consequence of the last proposition and the next.

\proclaim Proposition \sno{24}. Under the assumptions of Theorem
\sno13, $f=g$.

\Proof. The proof runs along the lines of Theorem \sno3.
 We use the inequality $f\le g$ and Lemma \sno21 without further
reference.
Let $\Pth_0\in\bPTh$ satisfy $f(\Pth_0) <
g(\Pth_0)$, and let $\Pth_1\in\bPTh\setminus\{\Pth_0\}$ be such that
 $$f(\Pth_0) = T_2 g(\Pth_0) = \max\{d_{01},d_{10}\} \ge d_{10}
  \quad\hbox{ and }\quad  g(\Pth_1) \ge d_{10}\,.$$
Then $f(\Pth_1) < d_{10}$
because otherwise, $g(\Pth_0) = T_2 f(\Pth_0) \le
\max\{d_{01},d_{10}\} = f(\Pth_0)$.

We also have $f(\Pth_1) <
g(\Pth_1)$ because otherwise we should have
 $f(\Pth_1) = g(\Pth_1)$, so that again $g(\Pth_0) = T_2f(\Pth_0)\le
 \max\{d_{01},d_{10}\}
 = f(\Pth_0)$.
Then, repeating all these steps, deduce the existence of some
$\Pth_2\in\bPTh \setminus \{\Pth_1\}$ such that
 $$f(\Pth_1) = T_2g(\Pth_1) = \max\{d_{12},d_{21}\} \ge d_{21} \quad\hbox{ and }
\quad g(\Pth_2) \ge d_{21},$$ and $f(\Pth_2) < d_{21}$ and $f(\Pth_2)
< g(\Pth_2)$. Combining these inequalities yields the relations
$$\max\{d_{01},d_{10}\} \ge d_{10} >  \max\{d_{12},d_{21}\} \ge d_{21}.$$
Since
$f(\Pth_0)=\max\{d_{01},d_{10}\}>f(\Pth_1)=\max\{d_{12},d_{21}\}$ we
get $\Pth_0\neq\Pth_1$.
Use induction to construct a
whole sequence $\Pth_0, \Pth_1, \ldots$ of points from $\bPTh$
satisfying
 $$\max\{d_{01},d_{10}\} \ge d_{10} >\max\{d_{12},d_{21}\} \ge d_{21}
  > \cdots > \max\{d_{n-1,n},d_{n-1,n}\} \ge d_{n,n-1} > \cdots$$
and $f(\Pth_n) = d_{n,n+1}$ for all $n\ge 0$.  In particular then,
$f(\Pth_n) > f(\Pth_{n+1})$, so the points $\P_n$ are all
different.  But this means that we have constructed a descending
chain of $\bPTh$ contrary to what is assumed in Theorem \sno{13}. Hence
there can be no $\Pth_0 \in\bPTh$ such that $f(\Pth_0) < g(\Pth_0)$.
\endproof


\def\secnum{5}
\sectitle{Stochastic models}%

In this section we prove the existence and uniqueness for Models 1 and 2
for a special class of point processes.
Let $\bf N$ denote the set of all countable sets
$\bPTh\subset\calX:=\bbbR^2\times[0,\pi)$ such that $\card(\bPTh\cap B\times[0,\pi))<\infty$
for all bounded sets $B\in\bbbR^2$. Any such $\bPTh$ is identified with a
(counting) measure $\card(\bPTh\cap\cdot)$. We equip $\bf N$ as usual with the
smallest $\sigma$-field $\bcal N$ making the mappings $\bPTh\mapsto\bPTh(C)$
measurable for all measurable $C\subset\calX$.
In this section and the next we consider a {\sl marked point process}
$\Psi$, that is a
random element in $\bf N$ defined on some abstract probability space
$(\Omega,\cal F,\bbbP)$ .
We make the following assumptions on $\Psi$.  Let $c$ be a finite positive
real number and $\bbbQ$ a probability measure on $[0,\pi)$. Then the $n\,$th
factorial moment measure $\alpha^{(n)}$ of $\Psi$
(see Daley and Vere-Jones (2008)) satisfies for each $n\in\bbbN$
$$
\alpha^{(n)}\big(\d(\P_1,\theta_1),\ldots,\d(\P_n,\theta_n)\big)
  \le c^n\d\P_1\cdots\,\d\P_n
 \,\bbbQ(\d \theta_1)\cdots \, \bbbQ(\d \theta_n), \seqno{1}
$$
where $\d\P$ denotes the differential of Lebesgue measure in $\bbbR^2$.
Assume also that the ground process $\Phi$, defined as the projection of
$\Psi$ on its first coordinate, is a simple point process. A stationary,
independently marked Poisson process with arbitrary mark distribution
satisfies \eqn1. So the mark distribution could for example be a sum of Dirac
measures as well as a diffuse measure. Moreover special classes of Cox and
Gibbs processes satisfy \eqn1. The details on this for the standard lilypond
model are stated in Daley and Last (2005) (= [D\&L]) and can be adapted to our
situation.

\proclaim Theorem \sno1. For $k=1,2$ and the marked point process $\Psi$ as
above, almost surely there exists a unique GMHS model of Type $k$.

We prove the theorem
by combining Propositions \sno2 and \sno3 with Theorem 4.3, and then
Proposition 5.4 with Theorem 4.13 for the cases $k=1$ and 2 respectively
(Theorems 4.3 and 4.13 from Section 4 show the growth-maximal property for
Models 1 and 2).  Consequently, Theorem \sno1 gives a precise meaning to
Models 1 and 2 described in the introduction.

\proclaim
Proposition \sno2. For the marked point process $\Psi$ as above, almost surely
there are no distinct pairs of points $\Pth_i, \Pth_j \in \Psi$ for which
$d_{ij} = d_{ji} <\infty$.

In other words, for a Poisson process Conditions D hold a.s.

\Proof.  The assertion can be proved as for Lemma 3.1 in 
[D\&L] showing a nonlattice property based on the factorial moment measure
condition on the point process.\endproof

\proclaim
Proposition \sno3.
For the marked point process $\Psi$ as above, almost surely there is no
descending chain of Type 1, i.e.\ there is no infinite sequence
$\Pth_0, \Pth_1, \Pth_2,\ldots$ of distinct points in $\Psi$ such that, with
$d_{ij} = d(\P_i, \P_{ij})$,
 $$\infty > d_{01} \ge d_{10} \ge \cdots \ge d_{n-1,n} \ge d_{n,n-1} \ge
  \cdots\,. \seqno{2}$$

\Proof.  We proceed as in Section 3.2 of [D\&L]. 
 Let $C$ be the set of all $\bPTh\in\bf N$ which
contain a descending chain and let $W_k:=[-k,k]^2$ be a square of
side length $2k$.  Furthermore let $B\subset \bbbR^2$ be a bounded
Borel set. For $s\le t$ and $B\in{\cal B}(\bbbR^2)$ let $C(n,s,t,B)$ be the set of all
$\bPTh\in\bf N$ whose projection on the first coordinate contains $n+1$
different points $\P_0, \P_1,\ldots,\P_n$ such that $\P_0\in B$ and
$t\ge d_{01} \ge d_{10} \ge \cdots \ge d_{n-1,n} \ge d_{n,n-1}\ge s$
and $C(s,t,B)$ the set of all $\bPTh\in\bf N$ whose projection (on
the first coordinate) contains an infinite series of points satisfying the
ordering condition at \eqn2    with $\P_0\in B$.  Moreover let
$C(s,t)$ be the set of all $\bPTh\in\bf N$ whose
projection contains an infinite series of points satisfying  the
ordering condition at \eqn2.    Clearly the sets $C(n,t,s,B)$ are
decreasing in $n$ and
$$ C(s,t,B)={\textstyle \bigcup_{n=1}^\infty} C(n,s,t,B),\quad
s\le t,\ B\in{\cal B}(\bbbR^2), $$
and $C(s,t,W_k)$ is increasing in $W_k$ with limit $C(s,t)$. It is
sufficient to show that there exists a sequence $\{t_i\}$ with
$\lim_{i\to\infty} t_i=\infty$ such that
$$
\lim_{n\to\infty} \bbbP\{\Psi\in C(n,t_i,t_{i+1},B)\} = 0
$$
for all bounded $B$ and
all $i$ because then, using the set identities given above,
$$
\bbbP\{\Psi\in C\}=\bbbP\big\{\Psi\in{\textstyle\bigcup_{i=1}^\infty\bigcup_{k=1}^\infty}
C(t_i,t_{i+1},W_k)\big\}\le \sum_{i=1}^\infty \bbbP\big\{\Psi\in{\textstyle
\bigcup_{k=1}^\infty} C(t_i,t_{i+1},W_k)\big\}=0
$$
Using assumption \eqn{1} on the factorial moment measures  of $\Psi$ we obtain as in
[D\&L] that $\bbbP\{\Psi \in C(n,s,t,B)\}$ is bounded by
 $$ c^n \!\int \!\cdots\!\int \bone\{\P_0\in B\}\,\bone\{t \ge
 d_{i-1,i} \ge d_{i,i-1}\ge s  \ \ (i=1,\ldots,n)\}
\,\d\P_0 \,\bbbQ(\d\theta_0)\ldots \d\P_n\, \bbbQ(\d\theta_n).  \seqno{3}
$$

Now let $D(n,s,t,B)$ be the set of all
$\bPTh\in\bf N$
whose projection contains $n+1$ different points $\P_0, \P_1,\ldots,\P_n$
such that $\P_0\in B$ and $ 
  t \ge d_{i-1,i} \ge d_{i,i-1}\ge s$ for $i=1,\ldots,n$.
Clearly $C(n,s,t,B)\subseteq D(n,s,t,B)$.
Therefore the expression at \eqn3 is bounded by
 $$\eqalign{
\! c^n \!\int \!\cdots\!\int \bone\{\P_0\in B\}
\,&\bone\{
t \ge d_{i-2,i-1} \ge d_{i-1,i-2}\ge s \ \ (i=1,\ldots,n-1)\} \  \cr
&\bone\{ t \ge d_{n-1,n} \ge d_{n,n-1}\ge s\}\,
\d\P_0 \,\bbbQ(\d\theta_0)\ldots \d\P_n\, \bbbQ(\d\theta_n).
}  \seqno{4}$$
This expression is bounded in turn by
 $$\eqalignno{
 c^n \!\int \!\cdots\!\int \bone\{\P_0\in B\}\,
&\bone\{
  t \ge d_{i-1,i} \ge d_{i,i-1}\ge s \ \ (i=1,\ldots,n-1)\} \cr
&\bone\{\P_n\in D(|\theta_n-\theta_{n-1}|,t-s) \}
   \ \d\P_0 \,\bbbQ(\d\theta_0)\ldots \d\P_n\, \bbbQ(\d\theta_n),
&\eqn{5}\cr}$$
where $D(\theta,x)$ is a diamond of side-length $x$ and inner angle $\theta$.
Now the volume of $D(\theta,l)$ is bounded by $x^2$, so
we can use Fubini's theorem to deduce that this expression is bounded by
 $$\eqalign{
\! 4\, (t-s)^2\,c^n\!\int \!\cdots\!\int \bone\{\P_0\in B\}
 &\bone\{ t \ge d_{i-1,i} \ge d_{i,i-1}\ge s \quad(i=1,\ldots,n-1)\}  \cr
& \qquad
\d\P_0 \,\bbbQ(\d \theta_0)\ldots \d\P_{n-1}\, \bbbQ(\d\theta_{n-1}). }
$$
Repeating this argument another $n-1$ times, the last expression
is bounded by
 $$4^n(t-s)^{2n}c^n \!\int \! \bone\{\P_0\in B\}\ \d\P_0 \,\bbbQ(\d \theta_0)
\,\le\, {[4c(t-s)^2]}^n\,\ell(B), \seqno{6}$$
so
 $\bbbP\{\Psi \in C(n,s,t,B)\}\le {[4c(t-s)^2]}^n \ell(B)$.
Choosing $t_0:=0$ and $t_{i+1}:=t_i + 1/\sqrt{5c}$ implies that the
right-hand side $\to 0$ as $n\to\infty$ geometrically fast, so the
proof is complete. \endproof

We now deduce Theorem \sno1 for the case $k=2$ by combining the next result
with Theorem 4.13 and get a precise meaning of Model 2.

\proclaim Proposition \sno4.   
For the random process based on the
marked point process $\Psi$ as above, almost surely
there is no descending chain of Type 2 in $\Psi$, i.e.\ there is no infinite
sequence
${\Pth_0}, {\Pth_1}, {\Pth_2}, \ldots$ of different points in $\Psi$
such that, with $d_{ij} = d(\P_i, \P_{ij})$,
 $$\infty >d_{10} \ge \max\{d_{1,2},d_{2,1}\}\ge d_{2,1}\ge \max\{d_{2,3},d_{3,2}\}
 \cdots  .$$
\smallskip

\Proof. The calculations are similar to those in the proof of Proposition
\sno3 except that we have to replace inequalities of the type
 $ t \ge d_{n-1,n} \ge d_{n,n-1}\ge s $
by
 $ t \ge \max\{d_{n-1,n},d_{n,n-1}\} \ge d_{n,n-1}\ge s $.
This leads to
 $\bbbP\{\Psi \in C(n,s,t,B)\}\le {[4c\,t(t-s)]}^n \ell(B)$.
Choosing $t_i:=\half\sqrt{i/c}$ yields
 $${[4c\, t_i(t_i-t_{i-1})]}^n
\ell(B)=\left(\frac{\sqrt{i}}{\sqrt{i}+\sqrt{i-1}}\right)^n
\ell(B)\le a^n \ell(B)$$
for some $a<1$ (and $a>\half$).
So $\lim_{n\to\infty} \bbbP\{\Psi \in C(n,s,t,B)\} = 0$ as before.
\endproof

\noindent{\bf Remark \sno5.} There are measurable mappings
$(\bPTh,\Pth)\mapsto R_k(\bPTh,\Pth)$ ($k=1,2$) from
${\bf N}\times\calX$ to $[0,\infty]$
such that the GMHS models of Type 1 and 2 in Proposition \sno4 are given
by $\{(\Pth,R_k(\Psi,\Pth)):\Pth\in\Psi\}$. These mappings
can be defined as the limit inferior of the recursions
in Section 4. We then have the useful translation invariance
 $$ R_k(\bPTh+P,\Pth+P)=R_k(\bPTh,\Pth), \quad P\in\bbbR^2, $$
where $\Pth+P$ denotes the translation of $\Pth$ in the first component
and $\bPTh+P:=\{\Pth+P:\Pth\in \bPTh\}$.
The measurability of $R_k$ has been implicitly assumed above.

\def\secnum{6}
\sectitle{Infinite clusters and percolation}%
In this section we fix a marked point process $\Psi$
with ground process $\Phi$.  Assume that $\Psi$
satisfies the factorial moment assumption (5.1), and that $\Psi$
is {\sl stationary}, i.e.\ for all $P\in\bbbR^2$ the distributions of $\Psi$
and $\Psi+P$ coincide, where $\Psi+P$ is the translation of $\Psi$ by $P$ in
the first component. The  {\sl intensity} of $\Psi$ (and of $\Phi$)
is defined by
$\lambda:=\bbbE \Phi([0,1]^2)$, which is the mean number of points of $\Phi$
in the unit square.  Assume $\Psi\ne \emptyset$ and $\lambda<\infty$.
We will show that a.s.\ there is no percolation in Model 2, i.e.\
there are no infinite clusters. Since Model 2
is akin to the lilypond model via contact between spherical
grains [DSS], we use the idea of a doublet; the earlier definition
can be rephrased here in our more formal language as follows.
Recall here the notation introduced in Remark 5.5.

\proclaim Definition \sno1. Two segment neighbours $\Pth,\Qth\in\Psi$
constitute a {\rm doublet} in Model 2  if\break
$R_2(\Psi,\Pth)=R_2(\Psi,\Qth)$.

 Thus, for a doublet pair $\{\Pth, {\rm Q}^\theta\}$,
$\Pth$ and ${\rm Q}^\theta$ are stopping segment neighbours of each other.


\proclaim Lemma \sno2.
Almost surely, in Model 2 every $\Pth\in \Psi$ has
at most one stopping segment neighbour.

\Proof.  When $\Pth_0\in\Psi$ has $\Pth_1\in\Psi$ as a stopping segment
neighbour, $R^{(2)}(\Psi,\Pth_0)=\max\{d_{01},d_{10}\} = m_{01}$.
For $\Pth_2$ also to
be a stopping segment neighbour of $\Pth_0$ then
$R^{(2)}(\Psi,\Pth_0)= m_{02}$.  By Conditions D, $m_{01}\ne m_{02}$, so we
have a contradiction.
\endproof

For the next result we need the following.
Define a graph on $\Psi\subset\calX$. Two nodes, i.e. two points of $\Psi$, share an edge if one is the
stopping segment neighbour of the other in the corresponding Model 2. Every
component of this graph is called a {\sl cluster}. This definition of a cluster
is consistent with our earlier definition in the introduction.  An immediate
consequence
of Lemma \sno2 is that {\sl every cluster has at most one doublet}.

\proclaim Lemma \sno3. Let $\Psi$ be a stationary marked point process
satisfying the factorial
moment measure condition. Then a.s.\
there does not exist any infinite cluster with a
doublet.

\Proof. The  statement is proved by adapting the argument in the
proof of [D\&L]'s Theorem 5.1.
\endproof

Here is the main result of this section.

\proclaim Theorem \sno4.  Let $\Psi$ be a stationary marked point process
satisfying the factorial
moment measure condition. Then a.s.\ there is no infinite cluster in Model 2.

\Proof.
Because of Lemma \sno3, it remains to show that there exists no infinite
cluster without a doublet, i.e.\  we have to show that a.s.\ there
does not exist an infinite sequence $\{\Pth_i\Colon i=0,1,\ldots\}$ such that
for every $i$, $\Pth_{i+1}$ is a stopping segment neighbour of $\Pth_i$
and $\{\Pth_i, \Pth_{i+1}\}$ is not a doublet.

Suppose on the contrary that such an infinite sequence exists.
Then invoking Conditions D when required and applying Definition 4.1(c) to
$\Pth_i$ with the stopping segment neighbour
$\Pth_{i+1}$ for $i=0,1,\ldots,$ we have 
 $$R_i = m_{i,i+1}\quad\hbox{ and }\quad d_{i+1,i} < R_{i+1} \le R_i\,,
\seqno{2}$$ which together imply that $R_i = d_{i,i+1} > d_{i+1,i}$
and hence that $d_{i-1,i} > d_{i,i+1}$.

 Let $B$ be a bounded Borel set.  Denote by $C'(n,s,t,B)$ the set
of all $\bPTh\in\bf  N$  whose projection contains $n+1$
different points $\P_0, \P_1,\ldots,\P_n$ such that $\P_0\in B$,
$t \ge d_{01} \ge d_{12} \ge \cdots \ge d_{n-1,n} \ge s$,
$R_i=d_{i,i+1}$ and $R_{i+1} \ge  d_{i+1,i}\,$.

Let $D'(n,s,t,B)$ be the set of all $\bPTh\in\bf N$
whose projection contains $n+1$ different points $\P_0,
\P_1,\ldots,\P_n$ such that $\P_0\in B$,
$t \ge d_{01} >d_{12} > \cdots > d_{n-1,n}\ge s$ and
$t\ge d_{i+1,i}$ for $1\le i\le n-1$.

Combining the last three conditions of the definition of $C'(n,s,t,B)$ we get
$$C'(n,s,t,B)\subseteq D'(n,s,t,B).$$
Analogously to the existence proof in Section 5, it is sufficient to
show that there exists a sequence $\{t_i\}$ with $\lim_{i\to\infty} t_i$
such that $\lim_{n\to\infty} \bbbP\{\Psi \in
D'(n,t_i,t_{i+1},B)\} = 0$ for all $B$ and all $i$.

As in Section 5, $\bbbP\{\Psi \in D'(n,t_i,t_{i+1},B)\}$ is bounded by
 $$\eqalignno{
\! c^n \!\int \!\cdots\!\int \bone\{\P_0\in B\}\,&\bone\{t \ge
d_{01} \ge d_{12} \ge \cdots \ge d_{n-1,n} \ge s \} \cr
&\bone\{t \ge d_{i+1,i},\ 0\le i\le n-1 \} \, \d\P_0
\,\bbbQ(\d\theta_0)\ldots \d\P_n\, \bbbQ(\d\theta_n). &\eqn3\cr}
$$
In turn this can be bounded by
$$\eqalignno{
\! c^n \!\int \!\cdots\!\int \bone\{\P_0\in B\}\,&\bone\{t \ge
d_{01} \ge d_{12} \ge \cdots \ge d_{n-2,n-1} \ge s \}  
 \bone\{t \ge d_{i+1,i},\ 0\le i\le n-2 \} \cr
&\bone\{t \ge d_{n-1,n}\ge s, \; t \ge d_{n,n-1}\}
\, \d\P_0 \,\bbbQ(\d\theta_0)\ldots \d\P_n\, \bbbQ(\d\theta_n). &\eqn4\cr}
$$
The integrand can be rewritten in terms of the maximum as in the
proof of Proposition 5.3 and we get the result in the same manner as
there.
\endproof

\vfill\eject

\def\secnum{7}
\sectitle{Finite clusters and discussion}%
The main concerns of this section are properties of a stationary lilypond
system of line-segments based on a stationary marked point process
$\Psi \ne \emptyset$ with intensity $\lambda$ and ground process $\Phi$
as in Sections 5 and 6 and for which
 the factorial moment assumption at (5.1) is satisfied.
Introduce a probability measure $\bbbP^0_\Phi$ (on
the underlying sample space) such that $\Psi$ has the {\sl Palm distribution}
$$
\bbbP^0_\Phi\{\Psi\in\cdot\}=
 \lambda^{-1}\bbbE \int_{[0,1]^2} \bone\{\Psi-\P\in\cdot\}\,\Phi(\d\P),
$$
where the shift $\Psi-\P$ of $\Psi$ has been defined in Remark 5.5 and integration with respect to $\Phi$ means
integration with respect to the associated counting measure.
This probability measure describes $\Psi$ as seen from a {\sl typical} point
of $\Phi$ (see Daley and Vere-Jones (2008) and Last (2010) for more detail
on Palm distributions). Note that  $\bbbP^0_\Phi\{0\in \Phi\}=1$.
If $\Psi$ is an independently marked stationary
Poisson process whose mark distribution
$\bbbQ$ has generic mark $R$, then the Slivnyak--Mecke theorem implies that
$\Psi\cup\{(0,R)\}$
has distribution $\bbbP^0_\Phi$ when $R$ is independent of $\Psi$.
Let
$\bbbE^0_\Phi$ denote the expectation operator with respect to $\bbbP^0_\Phi$.

For Model 1 we have not been able to resolve
whether or not the process of line-segments percolates in the Poisson case.
In Section 6 we showed the a.s.\ absence of percolation for Model 2.
This is not surprising because it resembles the standard
lilypond models of H\"aggstr\o m and Meester (1996) for which they showed there
is a.s.\ no percolation.
  We formulate our belief as follows.

\proclaim Conjecture \sno1.
In the Model 1 lilypond system of line-segments based on a stationary planar
Poisson process, there is a.s.\ no percolation.

This hypothesis was formulated on the basis of simulation work, and is
supported by its truth having been shown in the special 
case of lines oriented in just one of two directions by Christian Hirsch
(2013).
Evidence from simulations is based on
examining large numbers of realizations for finite systems of an
increasing number of points and recording the mean number of
points in the cluster to which the line-segment through the origin belongs.
In these we found no evidence of an increasing mean cluster size as might
be anticipated if a.s.\ an infinitely large cluster exists when there is an
infinite set of germs.

The conjecture can be cast as a random directed graph problem in which, for
each realization, the nodes are the points $\bP$ and each node $\P'$ say
has exactly one
outward-directed edge, namely to the node $\P''$ which is the centre of the
line-segment that stops the growth of the line-segment passing through $\P'$.
Resolving Conjecture \sno1 is the same as determining whether or not such a
graph can (with positive probability) have an infinitely large component.

Associate with each $\P_i$ of a realization of a system as in Conjecture \sno1
the vector $\bX_i := \P_{ij} \P_{jk}$, where $\Pth_i$ has $\Pth_j$ as its
stopping segment neighbour and $\Pth_j$ has $\Pth_k$ as its stopping
segment neighbour (in the notation of Algorithm 3.1, $j=J(i)$, $k=J(j)$).
Then tracing the successive `steps' $\{\bX_i\}$ within a cluster that has
no infinite line-segment, resembles tracing the steps of a random walk whose
mean step-length $\E(\bX_i) = 0$ (by rotational symmetry and the fact, from
Proposition \sno4 below, that $|\bX_i|$ has an exponentially bounded tail).
These steps are not independent (because of their construction), but they have
the property of successive steps ending in a cycle unless they are part of an
infinite cluster.  If we regard such `terminal' behaviour as indicating a
propensity for recurrence (as holds for a random walk in $\bbbR^2$ with no
drift), then this is further evidence to support Conjecture \sno1.

As in Section 2 call
a finite sequence $\Pth_{1},\dots,\Pth_n\in \Psi$ an {\sl$r$-cycle}
(in Model 1) if $\Pth_{i+1}$ is a stopping segment neighbour of
$\Pth_{i}$ for every $i=1,\dots,r$, where $\Pth_{r+1}:=\Pth_1$.
Clusters in Model 1 are as earlier.
It is easy to see that (almost surely) any finite cluster
has exactly one cycle while any infinite cluster
has at most one cycle. The next result is a first step
towards the proof of Conjecture \secnum.1.
Its proof is similar to the proof of Theorem 6.4.

\proclaim Proposition \sno2.
In Model 1 almost surely there is no infinite cluster with a cycle.

\Remark{7.3}.
We indicated in Section 1 that
there exists at least a third possible interpretation of
``growth-maximality'' with respect to hard-core models;
it is variously called a Gilbert tessellation or crack growth process
(see Schreiber and Soja (2011) for references).
In this model a line-segment stops growing only in the
direction in which it is blocked;
it is stopped in the other direction when it hits another line-segment.
Consequently, this model leads to a tessellation.
Schreiber and Soja (2011) prove  stabilization and a central limit theorem
for the Gilbert model. \endproof

While the two ends of a line-segment act ``independently''
of each other in this Gilbert model and that
is clearly not the case for our Models 1 and 2, one can prove the following
result      along the lines of Theorem 2.1 of Schreiber and Soja (2011).
Under $\bbbP^0_\Phi$ let $R^0$ denote the radius of the
(typical) line-segment centred at $0$.

\proclaim
Proposition \sno4.
Consider a stationary marked planar Poisson process with non-degenerate
mark distribution $\bbbQ^0$.
Then there are $\alpha,\beta>0$ such that
$$
\bbbP^0_\Phi\{R^0>t\}\le \alpha\exp(-\beta t^2),\quad t\ge 0.
$$

In the general case the (Palm) {\sl mark distribution} of $\Psi$ is the
probability measure $\bbbQ^0$ satisfying
$\bbbE[\Psi(\d(\P,\vartheta))]=\lambda\, \d \P\, \bbbQ^0(\d \vartheta)$.
We then have the following weak version of Proposition \sno4.

\proclaim
Proposition \sno5.
Let the process $\Psi$ of Proposition \sno4 be ergodic, and suppose that
$\bbbQ^0$ is diffuse.
Then a.s.\ there exists
no segment of infinite length, i.e.\   $\bbbP^0_\Phi\{R^0<\infty\}=1$.

\Proof. Let $\Psi^*:=\{\Pth\in\Phi: R(\Pth,\Psi)=\infty\} \subseteq \Psi$
denote the marked point process of
line-segments of infinite length,
where $R(\cdot,\cdot)$ refers to one of Models 1 and 2.  Observe that
$\{\Psi^*(\calX)=\infty\}$ is an invariant event so it has probability 0 or 1;
suppose for the sake of
contradiction that it has full probability.
Then there is a random direction $\vartheta\in[0,\pi)$ such that
all segments in $\Psi^*$ have this direction (the presence of a second
direction would contradict the hard-core property).  $\Psi$ is ergodic so
$\vartheta$ is non-random, and hence $\bbbQ^0$ has
an atom at $\vartheta$.  This is impossible for diffuse $\bbbQ^0$.
\endproof


For any $\P\in\Phi$ let $C(\P)\equiv C(\Psi,\P)$ denote the cluster
containing the line-segment centred at $\P$ and $\nu(\P)\equiv \nu(\Psi,\P)$
the number of neighbours of this line-segment.  For Model 1,
let $Z(\P)$ denote the unique
cycle${}\subseteq C(\P)$
when $\card\!\big(C(\P)\big) < \infty$, else set $Z(\P)=\emptyset$.
  For Model 2 let $D(\P)\subseteq C(\P)$
denote the doublet of $C(\P)$.  In developing certain mean value formulae in
the next two propositions we use
\vskip -19pt
 $$\eqalignno{
 \hbox{for Model 1,}\qquad \varpi_r &=\bbbP_\Phi^0\{O \in Z(O),
  \,\card Z(O) = r\} , &\eqn1\cr
 \hbox{and for Model 2},\;\qquad  \varpi&= \bbbP_\Phi^0\{O\in D(O)\},&\eqn2\cr
}$$
being the Palm probabilities that the line-segment through the origin $O$ is an
element of an $r$-cycle for Model 1 or an element of a doublet for Model 2.

\proclaim
Proposition \sno6. In Model 1, $\bbbE^0_\Phi \nu(O)=2$.
In Model 2,
$\bbbE^0_\Phi \nu(O)=2-\varpi$.

\Proof. For $\P,\Q\in\Phi$ let $\kappa(\P,\Q):=1$ if $\Q$ is a stopping
segment neighbour of $P$, $:=0$  otherwise. Let $B:=[0,1]^2$.
By the mass-transport principle (see e.g.\ Last (2010) equation (3.44)) we have
$$
\bbbE \int \!\int\bone_B(\P)\kappa(\P,\Q)\,\Phi(\d \Q)\,\Phi(\d \P)=
\bbbE \int \! \int\bone_B(\Q)\kappa(\P,\Q)\,\Phi(\d \P)\,\Phi(\d \Q).
\seqno3 $$
Because a.s.\ any line-segment has exactly one stopping neighbour the
left-hand side above equals the intensity $\lambda$.  For Model 1 the
right-hand side equals
 $$ \bbbE \int \!\int\bone_B(\Q)\,[\nu(\Q)-1]\,\Phi(\d \Q)
  =\lambda \,\bbbE^0_\Phi [\nu(O)-1],  $$
implying the first result.
The result for Model 2 comes from evaluating the right-hand side of \eqn1:
$$\eqalignno{
\bbbE &\int \bone_B(\Q)\,\bone\{\Q\in D(\Q)\}\nu(\Q)\,\Phi(\d \Q)+
\bbbE \int \bone_B(\Q)\,\bone\{\Q\notin D(\Q)\}\,[\nu(\Q)-1]\,\Phi(\d \Q) \cr
&=\lambda \,\bbbE^0_\Phi \nu(O)-\lambda \,\bbbP^0_\Phi\{O\notin D(O)\}
=\lambda \,\bbbE^0_\Phi \nu(O)-\lambda +\lambda p. &\endsqr
} $$

Because of the tree structure of any infinite cluster
and by analogy with a critical branching process,
Proposition \sno6 provides further evidence
supporting Conjecture \secnum.1.

Let $\Phi_c:=\{l(C(\P)):\P\in\Phi,\ \card(C(\P))<\infty\}$
denote the stationary point process of finite clusters, where
$l(A)$ denotes the lexicographic minimum of a finite set $A\subset\bbbR^2$;
let
$\lambda_c :=\break \bbbE\big[\card \Phi_c\big([0,1]^2\big)\big]$ denote
its intensity.
Because finite clusters are in one--one correspondence with cycles for Model 1
and doublets for Model 2, $\Phi_c$ can equally well be called a point process
of cycles or doublets.
Then
 $$ \mu:=\bbbE \int_{[0,1]^2} \card(C(\P))\,\Phi_c(\d \P) $$
can be interpreted as the mean size of the {\sl typical finite cluster}.

\proclaim
Proposition \sno7.
For Model 1 the mean size $\mu$ of the typical finite cluster is given by
 $$ \mu=\bbbP^0_\Phi\{\card C(O)<\infty\}
  \bigg(\sum^\infty_{r=3} \frac{\varpi_r}{r}\bigg)^{-1}
$$
and in Model 2 by $\mu=2/\varpi = (\varpi/2)^{-1}$, where $\varpi_r$ and
$\varpi$ are defined at \eqn1 and \eqn2.
\vfill\eject

\Proof. Consider Model 1. For $\P,\Q\in\Phi$ let
$\kappa(\P,\Q):=(\card Z(\P))^{-1}$ if
$\P\in Z(\P)$ and
$\Q=l(Z(P))$, $\kappa(\P,\Q):=0$ otherwise.
Then the left-hand side of \eqn3 equals
 $$ \lambda \,\bbbE^0_\Phi \bone\{O\in Z(O)\}\card(Z(O))^{-1}
   =\lambda\sum^\infty_{n=3} \frac{\varpi_r}{r}\,,
$$
and the right-hand side equals $\lambda_c=\bbbE \,\Phi_c\big([0,1]^2\big)$
because this intensity is the same as the intensity of the cycles.
Since $\mu$ is the quotient of the intensity of all points in finite clusters
and the intensity $\lambda_c$, and the first intensity is given by
$\lambda\bbbP^0_\Phi\{\card C(O)<\infty\}$, the first result follows.

The result for Model 2 follows by the same argument, as the intensity
of clusters is given by $\lambda \varpi/2$ and by Theorem 6.4
there are no infinite clusters.
\endproof

Last and Penrose (2012) established various properties for the standard
lilypond model in $\bbbR^d$, notably stabilizing properties,
a central limit theorem and frog percolation.  We believe that analogous
results should be available for both Models 1 and 2  for line-segments,
more easily for Model 2 because the techniques they used should continue to
be applicable.  A major task in adapting their proofs is to find an upper
bound on the length of a given line-segment as this may then be used to
replace the nearest-neighbour distance which they used as an upper bound on
the radius of a given hypersphere.

\bigskip\centerline{\sl References}\smallskip
\def\sc{\er}
\frenchspacing

\snh A{\sc NDRIENKO}, Yu.A., B{\sc RILLIANTOV}, N.V. and K{\sc RAPINSKY}, P.L.
(1994). Pattern formation by growing droplets: the touch-and-stop model
of growth. {\sl J. Statist. Phys. \bf75}, 507--523.

\snh D{\sc ALEY}, D.J., E{\sc BERT}, S. and S{\sc WIFT}, R.J. (2014). 
Size distributions in random triangles. {\sl J. Appl. Probab. \bf51A} 
(to appear).

\snh D{\sc ALEY}, D.J. and L{\sc AST}, G. (2005). \ Descending
chains, the lilypond model, and mutual-nearest-neighbour matching.
 {\sl Adv. Appl. Probab. \bf37}, 604--628.

\snh D{\sc ALEY}, D.J., M{\sc ALLOWS}, C.L. and S{\sc HEPP}, L.A.
(2000). \ A one-dimensional Poisson growth model with
non-overlapping intervals.\ {\sl Stoch. Proc. Appl. \bf90}, 223--241.

\snh D{\sc ALEY}, D.J., S{\sc TOYAN}, D. and S{\sc TOYAN}, H.
(1999). \ The volume fraction of a Poisson germ model with
maximally non-overlapping spherical grains.  {\sl Adv. Appl. Probab.
\bf31}, 610--624.

\snh D{\sc ALEY}, D.J. and V{\sc ERE}-J{\sc ONES}, D. (2003, 2008). \
{\sl An Introduction to the Theory of Point Processes} (Second Edition).
Volume I: Elementary Theory and Methods, and
Volume II: General Theory and Structure.
Springer, New York.

\snh
E{\sc BERT}, S. and L{\sc AST}, G.\ (2013). 
On a class of growth-maximal hard-core processes. arXiv:1303.2092,
submitted for publication.

\snh
G{\sc ILBERT}, E.N. (1967).  Random plane networks and needle-shaped
crystals.  In {\sl Applications of Undergraduate Mathematics in
Engineering},  B.\ Noble (Ed.), Macmillan, New York, Chap.\ 16.

\snh G{\sc RAY}, N.H., A{\sc NDERSON}, J.B., D{\sc EVINE}, J.D. and
K{\sc WASNIK}, J.M. (1976). \ Topological properties of random
crack networks. {\sl Math. Geology \bf 8}, 617--626.

\snh H{\sc \"AGGSTR\"OM}, O. and M{\sc EESTER}, R. (1996). \
Nearest neighbor and hard sphere models in continuum percolation.
{\sl Random Struct. Algorithms \bf 9}, 295--315.

\snh H{\sc EVELING}, M. and L{\sc AST}, G. (2006). \ Existence,
uniqueness, and algorithmic computation of general lilypond
systems.
{\sl Random Struct. Algorithms \bf 29}, 338--350.

\snh H{\sc IRSCH}, C. (2013).  Deterministic walks in random environment.
arXiv:1301.7279, submitted for publication.

\snh L{\sc AST}, G. (2010). \ Modern random measures: Palm theory and related
models. In {\sl New Perspectives in Stochastic Geometry\/} (eds. W.S. Kendall
and I. Molchanov), Oxford University Press, Oxford.

\snh L{\sc AST}, G. and P{\sc ENROSE}, M. (2013).
Percolation and limit theory for the Poisson lilypond model.
{\sl Random Struct. Algorithms \bf42}, 226--249.

\snh S{\sc CHREIBER}, T. and S{\sc OJA}, N. (2011).\ Limit
theory for planar Gilbert tesselations.
{\sl Probab. Math. Statist. \bf 31}, 149--160.


\snh S{\sc TIENEN}, J. (1982). \ {\sl Die Vergr\"oberung von Karbiden in reinen
Eisen-Kohlenstoff-Staehlen.}\/ Dissertation, RWTH Aachen.

\snh S{\sc TOYAN}, D.,  K{\sc ENDALL}, W.S. and M{\sc ECKE}, J. (1995).  \
{\sl Stochastic Geometry and its Applications}, 2nd Ed. John Wiley \& Sons,
Chichester.

\vfill\eject\end